\documentclass[review,onefignum,onetabnum]{siamart171218}

\usepackage{bm}
\usepackage{mathtools}
\usepackage{verbatim}
\usepackage{comment}
\usepackage{multirow} 
\usepackage{xspace,soul,color}
\usepackage{xcolor}
\soulregister{\cite}7 % 注册\cite命令
\soulregister{\citep}7 % 注册\citep命令
\soulregister{\citet}7 % 注册\citet命令
\soulregister{\ref}7 % 注册\ref命令

\usepackage{mathrsfs,amsmath,amssymb,bm}
\usepackage{lipsum}
\usepackage{amsfonts}
\usepackage{graphicx}
\usepackage{epstopdf}
\usepackage{makecell, rotating}
\usepackage{multirow}
\usepackage{algorithmic}
\usepackage{lipsum}
\usepackage{amsfonts}
\usepackage{graphicx}
\usepackage{epstopdf}
\usepackage{algorithmic}
\ifpdf
  \DeclareGraphicsExtensions{.eps,.pdf,.png,.jpg}
\else
  \DeclareGraphicsExtensions{.eps}
\fi

\newcommand{\bbR}{\mathbb{R}}

\newcommand{\bbZ}{\mathbb{Z}}

\newtheorem{thm}{Theorem}[section]

\newtheorem{remark}[thm]{Remark}

\newsiamremark{hypothesis}{Hypothesis}
\crefname{hypothesis}{Hypothesis}{Hypotheses}
\newsiamthm{claim}{Claim}

\headers{An efficient saddle search method for phase transitions }{Gang Cui, Kai Jiang, Tiejun Zhou }

\title{An efficient nullspace-preserving saddle search method for phase transitions involving translational invariance}

\author{Gang Cui\footnotemark[2]
        \and Kai Jiang\thanks{Hunan Key Laboratory for Computation and Simulation in Science and Engineering,
        Key Laboratory of Intelligent Computing and Information Processing of Ministry
        of Education, School of Mathematics and Computational Science, Xiangtan University, Xiangtan, Hunan,
        China, 411105.
        (Corresponding author. Kai Jiang, \email{kaijiang@xtu.edu.cn}).}
        \and Tiejun Zhou\footnotemark[2]}

\usepackage{amsopn}
\DeclareMathOperator{\diag}{diag}
\usepackage{amsfonts}

\makeatletter
\newcommand*{\addFileDependency}[1]{% argument=file name and extension
  \typeout{(#1)}% latexmk will find this if $recorder=0 (however, in that case, it will ignore #1 if it is a .aux or .pdf file etc and it exists! if it doesn't exist, it will appear in the list of dependents regardless)
  \@addtofilelist{#1}% if you want it to appear in \listfiles, not really necessary and latexmk doesn't use this
  \IfFileExists{#1}{}{\typeout{No file #1.}}% latexmk will find this message if #1 doesn't exist (yet)
}
\makeatother

\definecolor{myred}{rgb}{1,0.8,0.8}

\definecolor{mycyan}{rgb}{0.5,0.92,1.0}

\definecolor{mycyan}{rgb}{0.5,0.92,1.0}
\definecolor{mygreen}{rgb}{0.56,0.93,0.56}

\definecolor{myhl}{rgb}{1.0,0.98,0.56}
\sethlcolor{myhl}

\hypersetup{
	pdftitle={An efficient saddle search method for ordered phase transitions involving translational invariance},
	pdfauthor={Kai Jiang}
}

\begin{document}
	\maketitle
	
	% REQUIRED
	\begin{abstract}
%	The bottleneck of studying phase transitions is the barrier-crossing process composed of escaping from the basin of the local minimum and finding the saddle point. Breaking the bottleneck requires designing efficient algorithms relevant to the properties of concrete phase transition.
	In this work, we propose an efficient nullspace-preserving saddle search (NPSS) method for a class of phase transitions involving translational invariance, where the critical states are often degenerate. The NPSS method includes two stages, escaping from the basin and searching for the index-1 generalized saddle point. The NPSS method climbs upward from the generalized local minimum in segments to overcome the challenges of degeneracy. In each segment, an effective ascent direction is ensured by keeping this direction orthogonal to the nullspace of the initial state in this segment. This method can escape the basin quickly and converge to the transition states.
	We apply the NPSS method to the phase transitions between crystals, and between crystal and quasicrystal, based on the Landau-Brazovskii and Lifshitz-Petrich free energy functionals. Numerical results show a good performance of the NPSS method.
	%based on the relation between nullspace and translational invariance. Generally, the transition states in phase transitions correspond to index-1 generalized saddle points (GSPs), including degenerate and non-degenerate index-1 saddle points. For phase transitions involving translational invariance, degenerate local minima will lead to difficulties in escaping the basin. By keeping the ascent direction orthogonally constrained to the kernel space of the initial state, the NPSS method can obtain an ascent direction to efficiently escape from the basin and find the index-1 GSP. We apply the NPSS method to study nucleation and phase transition between crystals and between crystal and quasicrystal based on the Landau-Brazovskii and Lifshitz-Petrich models. Numerical results show that the NPSS method has a good performance in studying phase transitions with translational invariance. \see{Furthermore, through the study of MEP, we have identified an important critical point, the inflection point (IP), where symmetry breaking starts to occur during the phase transition. Before reaching the IP, the nullspaces on MEP have the $nullspace$-$preserving$ property.}
	\end{abstract}
	
	% REQUIRED
	
	% REQUIRED
	
	\section{Introduction}
	\label{sec:intr}

   Phase transition has been a topic of great interest for many years.
   For systems with energy landscape, phase transition is the most probable transition path, well-known as the minimum energy paths (MEPs), between the local minima of free energy. 
   Many computational approaches have been proposed to study phase transition, mainly including two classes, path-finding methods and surface-walking approaches. The former, such as the nudged elastic band approach\,\cite{henkelman2000improved,henkelman2000climbing,jonsson1998nudged}, the string method\,\cite{carilli2015truncation,e2002string,henkelman2000improved,lin2010numerical,ren2013climbing,samanta2013optimization}, and the equal-bond-length method\,\cite{cui2023finding} requires a suitable initialization connecting initial and final states. The latter starts from an initial metastable state and needs a proper ascent direction to search for the transition state. Representative methods contain the gentlest ascent dynamics\,\cite{e2011gentlest,gu2018simplified}, the dimer-type way\,\cite{henkelman1999dimer,zhang2012shrinking}, the step and slide algorithm\,\cite{miron2001step} and the high-index saddle dynamics (HiSD) method\,\cite{yin2021transition,yin2020construction,yin2019high,zhang2022error}. 
   
   The bottleneck of computing phase transitions involves escaping from the basin of a local minimum and locating the relevant saddle point due to the existence of energy barriers. To break this bottleneck, we need to design efficient algorithms by exploiting the properties of concrete phase transition.
   In present work, we pay attention to these phase transitions involving translational invariance. 
   Many important systems possess translational invariance including periodic crystal and quasiperiodic cases. The periodic systems have obviously translational symmetry. The quasiperiodic systems, including quasicrystals\,\cite{levine1984quasicrystals, shechtman1984metallic, tsai2000stable, wang2011origin, zeng2004supramolecular}, grained boundaries\,\cite{cao2021computing, jiang2022tilt,  sutton1995interfaces,van2002grain}, incommensurate structures\,\cite{bak1982commensurate,bruce1978theory, cao2018unconventional,sander1995incommensurate}, can be embedded into high-dimensional periodic systems, thus be viewed as being translational invariance in superspace. From the perspective of energy functional, these critical points (local minima, local maxima, saddle points) corresponding to ordered structures are usually $degenerate$. In other words, the Hessian of energy functional on the critical point has zero eigenvalues. The presence of nullspaces of Hessian makes escaping from the basin of local minimum difficult.

   In this work, we present a nullspace-preserving saddle search (NPSS) method for phase transitions with translational invariance. The NPSS method includes two stages, escaping from the basin and searching for the index-1 generalized saddle point (GSP). In Stage I, the NPSS method climbs upward from the generalized local minimum (GLM) in segments based on the changes of principal angles measuring the difference of nullspace. Specifically, within each segment, the ascent direction is chosen to be orthogonal to the nullspace of the initial state in the current segment. When the principal angle between nullspaces of the current and initial states exceeds a threshold, a new segment is initiated. These operations ensure the effectiveness of the ascent direction and avoid the costs of updating the nullspace at each step, enabling a quick escape from the basin. After escaping from the generalized quadratic region (GQR) of the GLM, the minimum eigenvalue of its Hessian becomes negative.  In Stage II, the NPSS method ascends along the ascent subspace $\mathcal{V}$ spanned by the eigenvector corresponding to the negative eigenvalue and descends along the orthogonal complement $\mathcal{V}^{\bot}$. Therefore, the NPSS method can effectively converge to the transition state.
    We demonstrate the power of the NPSS method for phase transitions with translational invariance by the Landau-Brazovskii (LB) and the Lifshitz-Petrich (LP) models. It is also noted that the NPSS can be applied to a range of models involving phase transitions with translational symmetry. 
    
    In this paper, mathematical derivations presented here are exemplified using the LP model. The article is structured as follows. In \Cref{sec:problem}, we provide an overview of the LB and LP models and employ the Fourier pseudo-spectral method for spatial discretization. In \Cref{sec:nullspace}, we analyze the properties of nullspaces.
   We present the NPSS method in \Cref{sec:method}. In \Cref{sec:result}, we apply the NPSS method to study the nucleation and phase transition between crystals, and between crystal and quasicrystal. We also compare with the HiSD method and demonstrate the superiority of the NPSS method in studying phase transitions with translational invariance.
    Finally, we discuss the relationship between symmetry-breaking and nullspaces on the MEP in \Cref{sec:discussion}. We find that the symmetry-breaking begins at the inflection point (IP), where the curvature changes sign on the MEP.
%	\newpage
	\section{Problem formulation}
	\label{sec:problem}
	Landau theory provides a powerful framework for studying the phases and phase transitions of ordered systems in both physical systems and materials\,\cite{cowley1980structural,provatas2011phase,toledano1987landau}. Specifically, let the order parameter function be $u(\bm{r})$, the Landau models can be expressed by a free energy functional
	\begin{align} 
		E(u; \Xi) = G(u; \Xi) + F(u; \Xi),
		\label{eq:landau_model}
	\end{align}
	where $\Xi$ represents physical parameters. $F(u)$ is the bulk energy with polynomial type formulation and $G(u)$ is the interaction energy that can contain higher-order differential operators to stabilize ordered structures.
	
	Two classes of Landau models are considered in this paper.
	The first one is the LB model, which can characterize the phases and phase transitions of periodic crystals\,\cite{brazovskii1975phase}. It has been applied in many different scientific fields, e.g., interfaces\,\cite{xu2017computing, Yao2022Transition} and polymeric materials\,\cite{shi1996theory,zhang2008efficient}. In particular, the energy functional of the LB model is
	\begin{align}
		E_{LB}(u) = \frac{1}{|\Omega|}\int_{\Omega}\left\{ \frac{1}{2}[(1 + \Delta)u]^2 +\frac{\tau}{2!}u^2 - \frac{\gamma}{3!}u^3 + \frac{1}{4!} u^4 \right\}\,d\bm{r}, \label{eq:LB}
	\end{align}
	where $u(\bm{r})$ is a real-valued function that measures the order of system. $\Omega$ is a bounded domain of the system with the volume $|\Omega|$. $\tau$ is the dimensionless reduced temperature, and $\gamma$ is the phenomenological coefficient. Compared with double-well bulk energy\,\cite{swift1977hydrodynamic}, the cubic term in the LB functional allows us to study the first-order phase transition.
	
	\begin{figure}[!ht]
		\centering
		\includegraphics[width=0.8\textwidth]{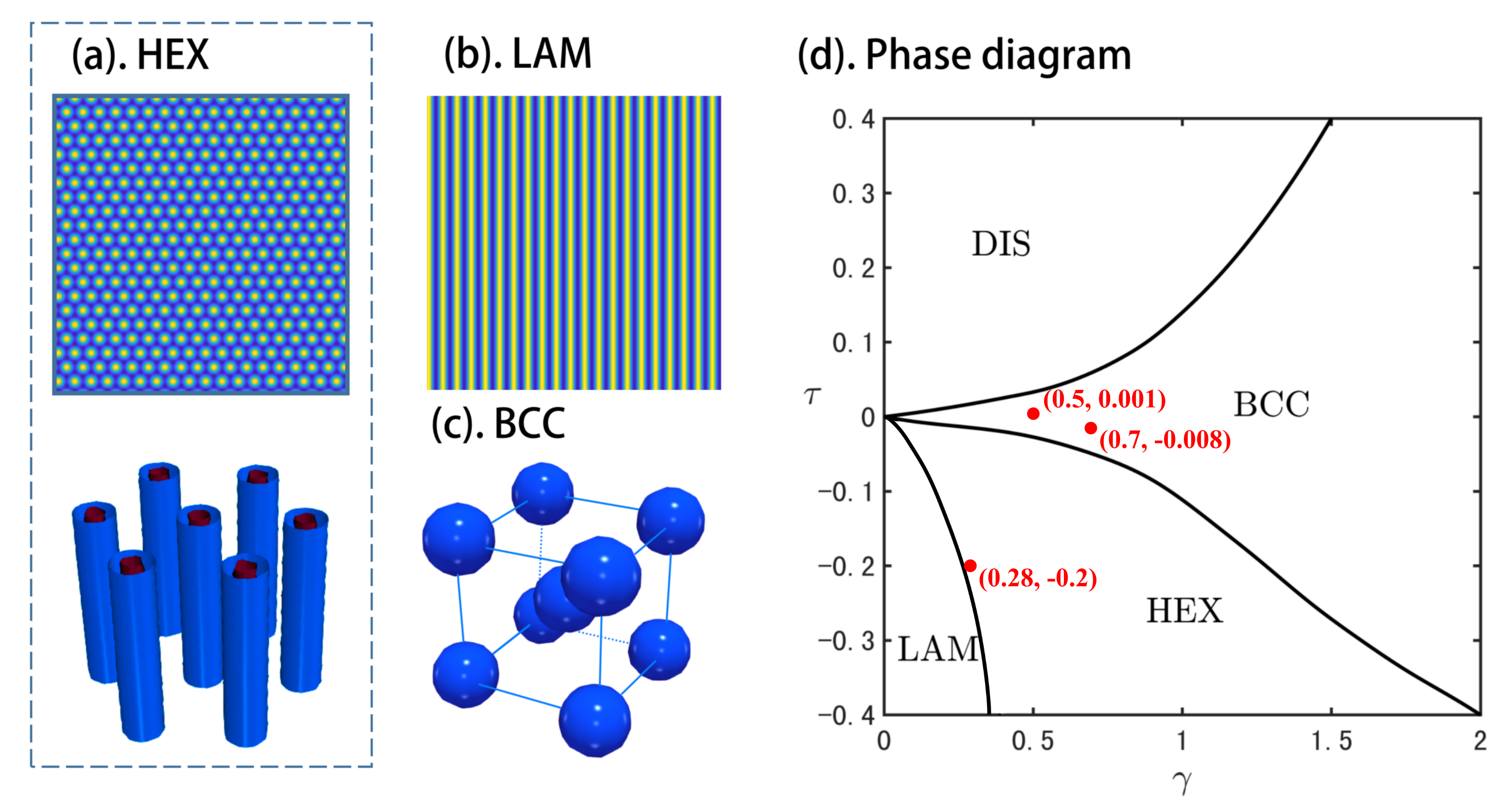}
		\caption{Stable ordered states and the phase diagram of the LB model.  Stable ordered structures:  (a). Hexagonal phase (HEX) in 2D and 3D; (b). Lamella phase (LAM); (c). Body-centered cubic spherical phase (BCC). DIS represents the disordered phase. Red dots represent the parameters $(\gamma, \tau)$ adopted in \Cref{sec:result}.}
		\label{fig:LB}
	\end{figure}
	The second one is the LP model, which can describe quasiperiodic structures, such as the bifrequency excited Faraday wave\,\cite{lifshitz1997theoretical}, and the stability of soft-matter quasicrystals\,\cite{jiang2015stability,lifshitz2007soft}. 
	\begin{align}
		E_{LP}(u) = \frac{1}{| \Omega |}\int_{\Omega}\left\{ \frac{1}{2}[(q_1^2 + \Delta)(q_2^2 + \Delta)u]^2 -\frac{\varepsilon}{2}u^2 - \frac{\alpha}{3}u^3 + \frac{1}{4} u^4 \right\}\, d\bm{r}, \label{eq:LP}
	\end{align}
	where $\varepsilon$ and $\alpha$ are phenomenological coefficients.  $q_1$ and $q_2$ are two characteristic length scales, which are necessary to stabilize the quasicrystals.
	\begin{figure}[!htb]
		\centering
		\includegraphics[width=0.8\textwidth]{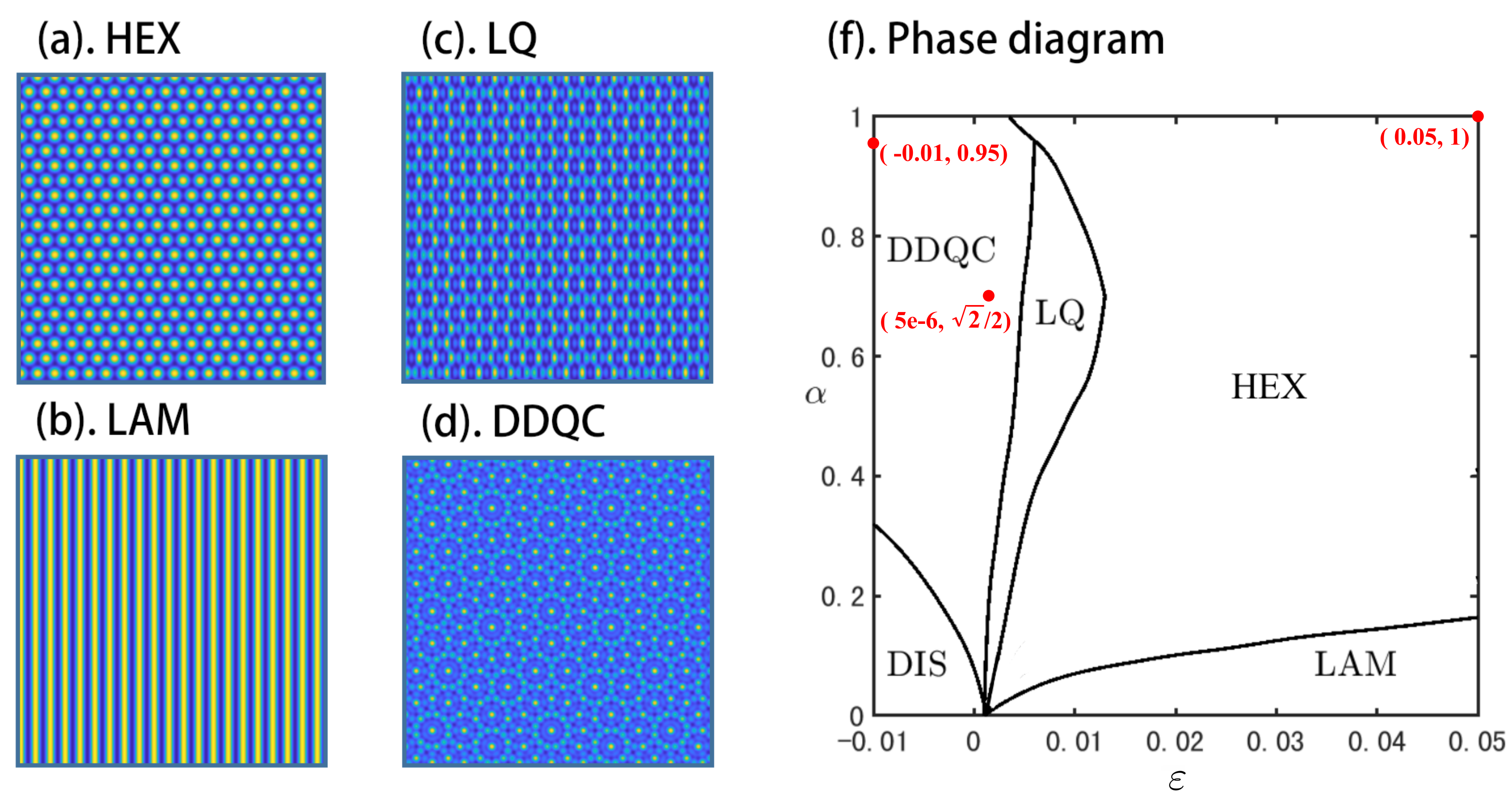}
		\caption{Stable ordered states and the phase diagram of the LP model with $q_1 = 1$, $q_2 = 2\cos(\pi/12)$. Stable ordered structures:  (a). HEX; (b). LAM; (c). Lamellar quasicrystal (LQ); (d). Dodecagonal quasicrystal (DDQC). Red dots represent the parameters $(\varepsilon, \alpha)$ adopted in \Cref{sec:result}.}
		\label{fig:LP}
	\end{figure}
	Furthermore, we impose the following mean zero condition of order parameter on the LB and LP systems, respectively, to ensure the mass conservation
	\begin{align}
		\frac{1}{|\Omega|}\int_{\Omega} u(\bm{r})\,d\bm{r}  = 0. \label{eq:constraint}
	\end{align}
	The equality constraint condition comes from the definition of the order parameter, which is the deviation from average density.
	
	Furthermore, \Cref{fig:LB} and \Cref{fig:LP} present the phase diagrams of the LB and LP models, respectively, only involving ordered structures considered in this paper. More sophisticated phase diagrams can refer to \cite{mcclenagan2019landau,yin2021transition}.
	
%	\newpage
	\subsection{Fourier pseudo-spectral method discretization} In this section, we introduce the Fourier pseudo-spectral method to discretize the LP energy functional. For an $n$-dimensional periodic order parameter $u(\bm{r})$, $\bm{r} \in \Omega:=\bbR^n / \bm{A}\bbZ^n $, where $\bm{A}=\left(\bm{a}_1, \bm{a}_2, \cdots, \bm{a}_n\right) \in \mathbb{R}^{n \times n}$ is the primitive Bravais lattice.
	The primitive reciprocal lattice $\bm{B} = \left(\bm{b}_1, \bm{b}_2, \cdots, \bm{b}_n\right) \in \mathbb{R}^{n \times n}$, satisfies the dual relation
	\begin{align}
		\bm{A}\bm{B}^{T}=2 \pi \bm{I}.
	\end{align}
	
	The order parameter $u(\bm{r})$ can be expanded as 
	\begin{align}
		u(\bm{r})= \sum_{\bm{k}\in \mathbb{Z}^n} \hat{u}\left(\bm{k}\right) e^{i (\bm{B} \bm{k})^T \bm{r}}, \quad \bm{r} \in \Omega,\label{eq:FourierSeries}
	\end{align}
	where the Fourier coefficient
	\begin{align}
		\hat{u}(\bm{k})=\frac{1}{|\Omega|} \int_{\Omega} u\left(\bm{r}\right) e^{-i (\bm{B} \bm{k})^T \bm{r}} d \bm{r}.
	\end{align}

	We define the discrete grid set as
	\begin{align}	\Omega_N=\Big\{(r_{1,j_1},\cdots,r_{n,j_n})=\bm{A}\big(j_1/N,\cdots,j_n/N \big)^T, j_i=0,\dots,N-1, i= 1,\cdots,n\Big\},\label{eq:discrete_grid}
	\end{align}
	where the number of elements of $\Omega_N$ is $M = N^n$. 
	Denote the grid periodic function space $\mathcal{G}_N$ = \big\{$f:\Omega_N \mapsto \mathbb{C}$, $f$ is periodic.\big\}. For any periodic grid functions $f,g \in \mathcal{G}_N$, we define the $\ell^2$-inner product as
	\begin{align}
		\langle f, g \rangle = \frac{1}{M}\sum_{\bm{r}_j \in \Omega_N} f(\bm{r}_j)\bar{g}(\bm{r}_j),
	\end{align}
	where $\bar{g}(\bm{r}_j)$ is the conjugate of $g(\bm{r}_j)$.
% 	And we define $\ell^2$ norm of $f$ as
    % \begin{align}
    %     \|f\|_{\ell^2} = \sqrt{\langle f, f \rangle}.
    % \end{align}
%     For any matrix $\bm{A} \in \mathbb{C}^{M \times M}$, we define its spectral norm as
%     \begin{align}
%         \|\bm{A}\|_2 = \sup_{\|f\|_{\ell^2} = 1} \|\bm{A}f\|_{\ell^2},
%     \end{align}
%     where $\bm{A}^{T}$ denotes the conjugate transpose of $\bm{A}$, and $\lambda_{\max}(\cdot)$ represents the largest eigenvalue of a matrix.
    The discrete reciprocal space is%MDPI: Please confirm if indentation needs to be added. Same as below.
	\begin{align}	\mathbb{K}_N^n=\{\bm{k}=\left(k_j\right)_{j=1}^n \in \mathbb{Z}^n:-N/2 \leq k_j<N/2\}.
	\end{align}
	For $\bm{k} \in \mathbb{Z}^n$ and $\bm{l} \in \mathbb{Z}^n$, we have the discrete orthogonality
	\begin{align}
		\langle e^{i (\bm{Bk})^T \bm{r_j}},e^{i (\bm{Bl})^T \bm{r_j}} \rangle = \left\{
		\begin{array}{cl}
			1, & \bm{k} = \bm{l}+N\bm{m}, ~\bm{m} \in \mathbb{Z}^n, \\
			0, & \mbox{otherwise}. 
		\end{array} \right.
		\label{eq:orthonormal}
	\end{align}
	Then the discrete Fourier coefficients of $u(\bm{r})$ in $\Omega_N$ can be represented as
	\begin{align}
		\hat{u}(\bm{k})= \langle u(\bm{r}_j), e^{i(\bm{Bk})^T\bm{r}_j}\rangle=\frac{1}{M} \sum_{\bm{r}_j\in \Omega_N} u\left(\bm{r}_j\right) e^{-i (\bm{B} \bm{k})^T \bm{r}_j}, \quad \bm{k}\in \mathbb{K}_N^n.
	\end{align}
	
	 Let $U = (u_1,u_2,\cdots,u_M) \in  \mathbb{R}^M$, where $u_j = u(\bm{r}_j), 1 \leq j \leq M$. Denote $\mathcal{F}_M \in  \mathbb{C}^{M\times M}$ as the discretized Fourier transformation matrix, then we have the discrete order parameter $\hat{U} = \mathcal{F}_M U$ in reciprocal space. The discretized energy function \cref{eq:LP} is
	\begin{align}
		E(\hat{U}) =G(\hat{U}) + F(\hat{U}),\label{eq:reciprocal_energy}
	\end{align}
	where $G$ and $F$ are the discretized interaction and bulk energies,
	\begin{align}
		G(\hat{U}) =& \frac{1}{2}\sum_{
			\bm{k}_1+\bm{k}_2=0}
		\bigl[ q_1^2 -(\bm{B k_1})^{\top}(\bm{Bk_1})\bigr]^2 \bigl[ q_2^2 -(\bm{B k_1})^{\top}(\bm{Bk_1})\bigr]^2  \hat{U}(\bm{k}_1)\hat{U}(\bm{k}_2), \label{eq:energies}\\
		F(\hat{U}) =& -\frac{\varepsilon}{2}\sum_{ 
			\bm{k}_1+\bm{k}_2= 0} \hat{U}(\bm{k}_1)\hat{U}(\bm{k}_2)- \frac{\alpha}{3} \sum_{\bm{k}_1+\bm{k}_2 + \bm{k}_3=0} 
		\hat{U}(\bm{k}_1)\hat{U}(\bm{k}_2)\hat{U}(\bm{k}_3)\notag\\
		&+ \frac{1}{4} \sum_{\bm{k}_1+\bm{k}_2 + \bm{k}_3 + \bm{k}_4 =0} 
		\hat{U}(\bm{k}_1)\hat{U}(\bm{k}_2)\hat{U}(\bm{k}_3)\hat{U}(\bm{k}_4),\notag
	\end{align}
	and $\bm{k}_j \in  \mathbb{K}_N^n$. The nonlinear terms $F(\hat{U})$ can be obtained by the FFT after being calculated in the $n$-dimensional physical space. Moreover, the mass conservation constraint \cref{eq:constraint} is discretized as
	\begin{align}
		e^{\top}_1 \hat{U} = 0, 
	\end{align}
	where $e_1 = (1, 0, \cdots , 0)^{\top} \in  \mathbb{R}^M$. Therefore, we obtain the following finite-dimensional
	minimization problem to obtain various ordered structures
	\begin{align}
		\min_{ \hat{U} \in \mathbb{C}^N} E(\hat{U}) = G(\hat{U}) + F(\hat{U}) \quad {\rm s.t.} \quad e^{\top}_1\hat{U} = 0.
	\end{align}
	According to \cref{eq:energies}, we have
	\begin{align}
		&\nabla G(\hat{U}) = D\hat{U}, \nabla F(\hat{U}) = \mathcal{F}_{M} \Lambda \mathcal{F}^{-1}_M \hat{U}, \label{eq:gradient}\\
		&\nabla^2 G(\hat{U}) = D, \nabla^2 F(\hat{U}) = \mathcal{F}_{M} \Lambda^{(')} \mathcal{F}^{-1}_M,\label{eq:hessian}
	\end{align}
	where $\mathcal{F}_M^{-1} \in \mathbb{C}^{M\times M}$ is the discretized Fourier inverse transform. $D$ is a diagonal matrix with non-negative entries $\bigl[ q_1^2 -(\bm{B k_1})^{\top}(\bm{Bk_1})\bigr]^2 \bigl[ q_2^2 -(\bm{B k_1})^{\top}(\bm{Bk_1})\bigr]^2$ and $\Lambda = \diag(-\varepsilon - \alpha U+ U^2)$, $\Lambda^{(')} = \diag(-\varepsilon - 2 \alpha U + 3U^2)$.

%	\newpage
	\section{Analysis of nullspace}
	\label{sec:nullspace}
 	 The nullspace of the Hessian on the MEP can affect the performance of saddle point search methods when studying phase transitions with translational invariance. In this section, we show the properties of nullspaces by theoretical analysis and numerical observations.
	
	Firstly, we rewrite the energy function \cref{eq:reciprocal_energy} in real space as 
	\begin{align}
		E(U) = G(U) + F(U), \label{eq:discrete_energy_function}
	\end{align}
	and the gradient and Hessian of $E(U)$ can be written as
	\begin{align}
		\nabla E(U) &= \mathcal{F}^{-1}_M(\nabla G(\hat{U}) +\nabla F(\hat{U})) = \mathcal{F}^{-1}_M D \hat{U}+\Lambda \mathcal{F}^{-1}_M \hat{U} = \mathcal{F}^{-1}_M D \mathcal{F}_M U +\Lambda U,\label{eq:discrete_gradient}\\
		H(U) &=\nabla^2 E(U) = \mathcal{F}^{-1}_M D \mathcal{F}_M +\Lambda^{(')}.\label{eq:discrete_hessian}
	\end{align}
 
	Now we discuss the critical points of degenerate and non-degenerate cases. For an energy function $E$ represented by \cref{eq:discrete_energy_function}, a critical point of $E$ satisfies $\nabla E(U) = 0$. To distinguish different critical points mathematically, we introduce the Morse index\,\cite{stein1963morse}. The Morse index of a non-degenerate critical point $U$ is equal to the number of negative eigenvalues of $H(U)$. For degenerate problems, we can extend the definition of the Morse index to depict critical points.
	\begin{definition}\label{generalized_critical_point}
		(1). Generalized Morse index. For a degenerate or non-degenerate critical point, the generalized Morse index is the number of negative eigenvalues of its Hessian. \\
		(2). Generalized critical point. A critical point is called the index-$k$ GSP if its generalized Morse index is $k (k > 0)$. If $k =1$, it is an index-1 GSP, characterized as the transition state. If $k=0$, it is a GLM.
	\end{definition} 
	Notably, the Morse index is a special case of the generalized Morse index. Similarly, we define the GQR of the critical points,
	
	\begin{definition}\label{generalized_quadratic_region}
	A region around a generalized critical point $U^{*}\in\mathbb{R}^M$ is the GQR, if the Hessian $H(U)$ has the same number of negative eigenvalues as $H(U^{*})$ for any $U\in\mathbb{R}^M$ in this region.
	\end{definition} 
	The quadratic region of a non-degenerate critical point has been mentioned in\,\cite{nagyfalusi2022magnetic}. Moreover, we define unstable,
	stable subspaces and nullspace of $H(U)$ as
	\begin{align}
		\mathcal{W}^u(U) &= {\rm span}_{\mathbb{R}}\{\bm{v}^u_1,\cdots,\bm{v}^u_{l_u}\},\\
		\mathcal{W}^s(U) &=  {\rm span}_{\mathbb{R}}\{\bm{v}^s_{1},\cdots,\bm{v}^s_{l_s}\},\\
		\mathcal{W}^k(U) &= {\rm span}_{\mathbb{R}}\{\bm{v}^k_{1},\cdots,\bm{v}^k_{l_k}\},
	\end{align}
	where $\{\bm{v}^u_1,\cdots,\bm{v}^u_{l_u}\}$, $\{\bm{v}^s_{1},\cdots,\bm{v}^s_{l_s}\}$ and $\{\bm{v}^k_{1},\cdots,\bm{v}^k_{l_k}\}$ $\subset \mathbb{R}^M$ are eigenvectors of $H(U)$
	corresponding to negative, positive, and zero eigenvalues, respectively. $l_u + l_s + l_k = M$.
	According to the primary decomposition theorem \cite{mh1974differential}, $\mathbb{R}^M$ can be decomposed as
	\begin{align}
		\mathbb{R}^M = \mathcal{W}^u(U) \oplus \mathcal{W}^s(U) \oplus \mathcal{W}^k(U), \label{eq:direct_sum}
	\end{align}
	where the dimensions $l_u$, $l_s$ and $l_k$ of subspaces $\mathcal{W}^u(U)$, $\mathcal{W}^s(U)$ and $\mathcal{W}^k(U)$ depend on the location of $U$ on the energy surface. 

    Define the norm $\|\bm{u}\|_2 = \sqrt{\sum_{i=1}^M u_i^2}$ for $\bm{u}=(u_1,\cdots,u_M ) \in \mathbb{R}^M$. To measure the difference between two subspaces, we give the definition of principal angles\,\cite{ake1973numerical,miao1992principal}.
	\begin{definition}\label{def:principal-angles}
	Let $\mathcal{W}$ and $\widehat{\mathcal{W}}$ be subspaces in $\mathbb{R}^M$, with {\rm dim} $\mathcal{W} = d \leq$ {\rm dim} $\widehat{\mathcal{W}} = m$. Then the principal angles between $\mathcal{W}$ and $\widehat{\mathcal{W}}$,
	$0 \leq \theta_1 \leq \theta_2 \leq \cdots \leq \theta_d \leq \pi/2$,
    are defined by
	\begin{align}
    \cos \theta_i := \frac{\bm{u}_i^{T} \bm{v}_i}{\|\bm{u}_i\|_{2}\|\bm{v}_i\|_{2}}
     = \max \left\{ \frac{\bm{u}^{T} \bm{v} }{\|\bm{u}\|_{2}\|\bm{v}\|_{2}} :\begin{array}{cc}
          \bm{u} \in \mathcal{W}, & \bm{u} \bot \bm{u}_k,\\
          \bm{v} \in \widehat{\mathcal{W}}, & \bm{v} \bot \bm{v}_k,
     \end{array} k = 1,\cdots,i-1\right\},\notag
	\end{align}
	where
	\begin{align}
    (\bm{u}_i, \bm{v}_i) \in \mathcal{W} \times \widehat{\mathcal{W}}, \quad i =  1,\cdots, d,
	\end{align}
	are the corresponding d pairs of principal vectors.
\end{definition}
    We also denote  
    $\sin \Theta(\mathcal{W}, \widehat{\mathcal{W}})$ as 
    \begin{align}\label{eq:sin_theta}
        \sin \Theta(\mathcal{W}, \widehat{\mathcal{W}}) &= \operatorname{diag}(\sin \theta_1, \sin \theta_2, \cdots, \sin \theta_d),
    \end{align}
    which is a $d \times d$ diagonal matrix\,\cite{ake1973numerical}. 
    Since $\sin \Theta(\mathcal{W}, \widehat{\mathcal{W}})$ is diagonal, its spectral norm $\|\sin \Theta(\mathcal{W}, \widehat{\mathcal{W}})\|_2$ is equal to the maximum of its diagonal elements, i.e.,
    \begin{align}\label{eq:norm_sin_theta}
        \|\sin \Theta(\mathcal{W}, \widehat{\mathcal{W}})\|_2 =  \sin \theta_d.
    \end{align}
    This norm can measure the difference between subspaces $\mathcal{W}$ and $\widehat{\mathcal{W}}$. 
    
	\subsection{Translational invariance and nullspace}\label{sec:periodicity}
	Critical points involved in phase transitions with translational invariance generally correspond to periodic or quasi-periodic phases. Although quasicrystals lack translational invariance in physical space, they can be embedded in higher-dimensional periodic systems. Thus we can analyze the properties of these phases uniformly in suitable periodic systems. For a non-zero periodic critical point, we discuss the relationship between nullspace's dimension of its Hessian and translational invariance in this subsection\,\cite{bao2024convergence}.
	
	\begin{proposition}
		Let $U$ be a $n$-dimensional non-zero periodic critical point of $E$ \cref{eq:discrete_energy_function}, i.e.,  $U (U \neq 0)$ is periodic with $\nabla E(U)= 0$, then $\dim \mathcal{W}^k(U) \geq n$. \label{prop:kernel_space}
	\end{proposition}
	
	\begin{proof}
		According to \cref{eq:discrete_gradient}, the critical point $U$ ($U \neq 0$) satisfies
		\begin{align}
			\nabla E(U)= \mathcal{F}^{-1}_M D \mathcal{F}_M U +\Lambda U = 0. \label{eq:stable_condition}
		\end{align}
		Define the discrete derivative operator as \begin{align}
		\mathcal{F}^{-1}_MK_j\mathcal{F}_M, \quad j=1,2,\cdots,n, \label{eq:derivative_operator}
		\end{align}
		where $K_j$ is a diagonal matrix, the diagonal elements relate to $i k_j$, where $k_j$ is the $j$-th component of spectral point $\bm{k}$ and $i$ is the imaginary unit. 
		
		Acting operator \cref{eq:derivative_operator} on \cref{eq:stable_condition} leads to 
		\begin{align}
		\mathcal{F}^{-1}_M K_j D\mathcal{F}_MU + \mathcal{F}^{-1}_M K_j\mathcal{F}_M\Lambda U =		\mathcal{F}^{-1}_M K_j D\hat{U} + \mathcal{F}^{-1}_M K_j\widehat{\Lambda U} =0.\label{eq:Derivative_gradient}
		\end{align}
		Let 
		\begin{align}
			\bm{v}_j = \mathcal{F}^{-1}_MK_j\hat{U} \in \mathbb{R}^M, \quad j= 1, 2, \cdots, n. \label{given_v}
		\end{align}
		By acting the Hessian \cref{eq:discrete_hessian} on $\bm{v}_j$,
		\begin{align}\label{eq:discrete_hessian_v}
        & (\mathcal{F}^{-1}_M D\mathcal{F}_M + \Lambda^{(')})\bm{v}_j
        =(\mathcal{F}^{-1}_M D\mathcal{F}_M + \Lambda^{(')})\mathcal{F}^{-1}_M K_j\hat{U}  \\
		&=  \mathcal{F}^{-1}_M D\mathcal{F}_M\mathcal{F}^{-1}_M K_j\hat{U} + \Lambda^{(')}\mathcal{F}^{-1}_M K_j\hat{U} \notag \\
		&=  \mathcal{F}^{-1}_M D K_j \hat{U} + \Lambda^{(')}\mathcal{F}^{-1}_M K_j\hat{U} \notag \\
		&\mathop{=}\limits^{(*)}\mathcal{F}^{-1}_M K_j D \hat{U} + \Lambda^{(')}\mathcal{F}^{-1}_M K_j\hat{U},\notag 
		\end{align}
		where $(*)$ holds because $D$ and $K_j$ are diagonal, $D K_j = K_j D$.
		
		Denote $P = (p_s)_{s=1}^M = \mathcal{F}^{-1}_M K_j\widehat{\Lambda U}$ and $Q = (q_s)_{s=1}^M  = \Lambda^{(')}\mathcal{F}^{-1}_M K_j\hat{U} $. Let $\bm{k},~ \bm{k}_1,~ \bm{k}_2,~ \bm{k}_3\in \mathbb{K}_N^n$, and $k_j,~ k_{1,j},~ k_{2,j},~ k_{3,j}$ be the $j$-th component of them, respectively. $\bm{r}_s = (r_1,\cdots,r_n) \in \Omega_N$ is the coordinate of the $s$-th element in $U$. Denote $p_s$ and $q_s$ as the $s$-th element of $P$ and $Q$. For $s =1,2,\cdots,M$, $p_s$ and $q_s$ can be written as
	\begin{align}
	p_s = & -\varepsilon\sum_{\bm{k} \in \mathbb{K}_N^n} ik_j \hat{u}(\bm{k})e^{i \bm{k} \bm{r}_s} - \alpha \sum_{\bm{k} \in \mathbb{K}_N^n} i k_j\sum_{\bm{k}_1 +\bm{k}_2 = \bm{k}}\hat{u}(\bm{k}_1)\hat{u}(\bm{k}_2)e^{i \bm{k} \bm{r}_s} \label{eq:p_s}\\
	&+\sum_{\bm{k} \in \mathbb{K}_N^n} ik_j\sum_{\bm{k}_1 +\bm{k}_2 +\bm{k}_3 = \bm{k}}\hat{u}(\bm{k}_1)\hat{u}(\bm{k}_2)\hat{u}(\bm{k}_3)e^{i \bm{k} \bm{r}_s}\notag \\
	= & -\varepsilon\sum_{\bm{k} \in \mathbb{K}_N^n} ik_j \hat{u}(\bm{k})e^{i \bm{k} \bm{r}_s} - \alpha \sum_{\bm{k} \in \mathbb{K}_N^n} i(k_{1,j}+k_{2,j})\sum_{\bm{k}_1 +\bm{k}_2 = \bm{k}}\hat{u}(\bm{k}_1)\hat{u}(\bm{k}_2)e^{i \bm{k} \bm{r}_s} \notag\\
	&+\sum_{\bm{k} \in \mathbb{K}_N^n}i(k_{1,j}+k_{2,j}+k_{3,j})\sum_{\bm{k}_1 +\bm{k}_2 +\bm{k}_3 = \bm{k}}\hat{u}(\bm{k}_1)\hat{u}(\bm{k}_2)\hat{u}(\bm{k}_3)e^{i \bm{k} \bm{r}_s}\notag \\
	\mathop{=}\limits^{(\#)} & -\varepsilon\sum_{\bm{k} \in \mathbb{K}_N^n} i k_j \hat{u}(\bm{k})e^{i \bm{k} \bm{r}_s} - 2 \alpha \sum_{\bm{k} \in \mathbb{K}_N^n} \sum_{\bm{k}_1 +\bm{k}_2 = \bm{k}}  i k_{1,j} \hat{u}(\bm{k}_1)\hat{u}(\bm{k}_2)e^{i \bm{k} \bm{r}_s} \notag\\
	&+ 3 \sum_{\bm{k} \in \mathbb{K}_N^n}\sum_{\bm{k}_1 +\bm{k}_2 +\bm{k}_3 = \bm{k}} i k_{1,j} \hat{u}(\bm{k}_1)\hat{u}(\bm{k}_2)\hat{u}(\bm{k}_3)e^{i \bm{k} \bm{r}_s},\notag\\
	\notag\\
	q_s =& \bigl[-\varepsilon - 2\alpha \sum_{\bm{k} \in \mathbb{K}_N^n}\hat{u}(\bm{k})e^{i \bm{k} \bm{r}_s} \label{eq:q_s}\\ &+3\sum_{\bm{k} \in \mathbb{K}_N^n}\sum_{\bm{k}_2 +\bm{k}_3 = \bm{k}}\hat{u}(\bm{k}_2)\hat{u}(\bm{k}_3)e^{i \bm{k} \bm{r}_s}\bigr]\sum_{\bm{k} \in \mathbb{K}_N^n} i k_j\hat{u}(\bm{k})e^{i \bm{k} \bm{r}_s}\notag\\
	= &-\varepsilon\sum_{\bm{k} \in \mathbb{K}_N^n}i k_j\hat{u}(\bm{k})e^{i \bm{k} \bm{r}_s} - 2\alpha \sum_{\bm{k} \in \mathbb{K}_N^n}\sum_{\bm{k}_1 +\bm{k}_2 = \bm{k}}ik_{1,j}\hat{u}(\bm{k}_1)\hat{u}(\bm{k}_2)e^{i \bm{k} \bm{r}_s}\notag\\
	& + 3\sum_{\bm{k} \in \mathbb{K}_N^n}\sum_{\bm{k}_1 +\bm{k}_2 +\bm{k}_3 = \bm{k}}ik_{1,j}\hat{u}(\bm{k}_1)\hat{u}(\bm{k}_2)\hat{u}(\bm{k}_3)e^{i \bm{k} \bm{r}_s},\notag
\end{align}
	where $(\#)$ holds because $k_{1,j}$, $k_{2,j}$ and $k_{3,j}$ are symmetric.
	
	From \cref{eq:p_s} and \cref{eq:q_s}, we have $p_s = q_s$.
	Therefore, 
	\begin{align}
	\mathcal{F}^{-1}_M K_j\widehat{\Lambda U} = \Lambda^{(')}\mathcal{F}^{-1}_M K_j\hat{U}.
	\end{align}
	Combined with \cref{eq:Derivative_gradient} and \cref{eq:discrete_hessian_v}, we have
        \begin{align}
        (\mathcal{F}^{-1}_M D\mathcal{F}_M + \Lambda^{(')})\bm{v}_j = 0,
        \end{align}
	    then we can conclude that $\bm{v}_j \in$ $\mathcal{W}^k(U), j=1,2,\cdots,n$. 
	    
		Next, we consider the dimension of ${\rm span} \{\bm{v}_1, \cdots, \bm{v}_n \}$. Denote $A\circ B =[a_{ij}b_{ij}]$ as the Hadamard product, where $A=[a_{ij}]$ and $B=[a_{ij}]$ are two matrices with the same order. Then
		\begin{align}
        \mathcal{F}_M\bm{v}_j = K_j\hat{U} = i  \tilde{\bm{k}}_j \circ \hat{U}, \quad j= 1,2,\cdots,n,
        \end{align}
        where $\tilde{\bm{k}}_j(s) =K_j(s,s)$. Define
		$\bm{v} = \sum^n_{j=1} \alpha_j \bm{v}_j$ with $\bm \alpha = (\alpha_1 , \cdots , \alpha_n )^{\top} \in \mathbb{R}^n$. $\mathcal{F}_M \bm{v}$ are give by
		\begin{align}
			\mathcal{F}_M \bm{v} = (\mathcal{F}_M\bm{v}_1, \cdots,\mathcal{F}_M\bm{v}_n)\bm{\alpha}
			 = i  (\tilde{\bm{k}}_1\circ \hat{U},\cdots, \tilde{\bm{k}}_n \circ \hat{U})\bm{\alpha}
			 = i  \hat{U}\circ(\tilde{\bm{k}}_1,\cdots,\tilde{\bm{k}}_n)\bm{\alpha}.
		\end{align}
		
		Denote $\bm{k} = (k_1,\cdots,k_n) \in \mathbb{K}_N^n$ as any row of $(\tilde{\bm{k}}_1,\cdots,\tilde{\bm{k}}_n)$, and $\mathcal{K}(U) := \{\bm k \in \mathbb{K}^n_N : \hat{U}(\bm k) \neq 0\}$, we have
		\begin{align}
			\bm{v} = 0 \Longleftrightarrow \mathcal{F}_M \bm{v}=0 {\rm \enspace for \, all} \enspace  \bm k \in \mathbb{K}^n_N \Longleftrightarrow
			\bm{k}\cdot \bm{\alpha} = 0 {\rm \enspace for \, all} \enspace \bm k \in \mathcal{K}(U),
		\end{align}
		implying that the dimension of ${\rm span} \{\bm{v}_1, \bm{v}_2, \cdots, \bm{v}_n \}$ is exactly the dimension
		of ${\rm span } \mathcal{K}(U)$. If $U$ is periodic on $\Omega_N$, all reciprocal lattice vectors of $U$ can be linearly expressed by $n$  basis\,\cite{jiang2014numerical}. Then we know that {\rm dim} ${\rm span }  \mathcal{K}(U) = n$, hence {\rm dim} $\mathcal{W}^k(U) \geq n$. The proof is completed.
	\end{proof}
	
	\begin{remark}
	\Cref{prop:kernel_space} reveals the relation between the dimension of nullspace and the periodicity of a non-zero critical point. To verify this conclusion, \cref{tab:eigenvalues} presents seven smallest eigenvalues of Hessian corresponding to three ordered structures.
	\renewcommand{\arraystretch}{1.5}
		\begin{table}[H]
		\centering
		\footnotesize{
		\caption{Seven numerically smallest eigenvalues of $H(U)$ corresponding to stable ordered structures in the LP model. Model parameters $\varepsilon = 0.05$, $\alpha = 1$, computational domain $\Omega = [0, 112 \times 2\pi)^2$, and spatial discretization points $N = 1024$ in each dimension.}}
		\begin{tabular}{|c|c|c|c|c|c|c|c|}  
		\hline  
		HEX  &-6.53e-17 &-1.26e-17  &5.57e-4  &5.567e-4 &5.58e-4  &5.58e-4   &1.11e-3   \\
		\hline  
		LQ   & -9.50e-16 &-1.28e-17  &1.82e-17   &2.57e-4   &2.57e-4 &5.80e-4  &5.80e-4          \\  
		\hline
        DDQC & -5.42e-17   &-5.99e-17    &-1.44e-17    &-7.67e-18   & 1.69e-3    &1.69e-3    &1.69e-3         \\
		\hline  
		\end{tabular}
		\label{tab:eigenvalues}
	\end{table}
	In this case, HEX is a 2-dimensional periodic structure, and LQ, DDQC can be embedded into 3- and 4-dimensional periodic systems, respectively. From \cref{tab:eigenvalues}, HEX, LQ, and DDQC have 2, 3, and 4 zero eigenvalues, respectively. It can be observed that the dimension of the embedded periodic system is exactly equal to the dimension of the corresponding nullspace. The same phenomenon has also been observed in LP and LB models for other periodic structures.
	It seems to have a more stronger conclusion, {\rm dim} $\mathcal{W}^k(U) = n$, although the definitive proof remains open. 
	\end{remark}

	\subsection{Impact of nullspace on escaping from basin} \label{sec:stability_energy}
	In ordered phase transition, our target is to locate an index-1 GSP from a GLM. However, due to the degeneracy of GLMs, the presence of nullspace hinders the escape from the basin of attraction. Therefore we analyze the effect of nullspace on escaping from the basin in this subsection. 
	
	Define $\|\cdot\|_{\ell^2} = \sqrt{\langle \cdot, \cdot\rangle}$. For a non-zero periodic GLM $\tilde{U}$, the subspaces of its Hessian include $\mathcal{W}^s(\tilde{U})$ and $\mathcal{W}^k(\tilde{U})$. The energy remains stable when perturbed in $\mathcal{W}^k(\tilde{U})$ at $\tilde{U}$.
	
	\begin{proposition}
		Let $\tilde{U}$ be a non-zero periodic GLM of $E(U)$. For a given $\epsilon > 0$,  any $\delta U \in \mathcal{W}^k(\tilde{U})$ satisfies $\|\delta U\|_{\ell^2} < \epsilon $, we have $|E(\tilde{U} + \delta U) - E(\tilde{U})| < o(\epsilon^2)$. \label{prop:stability_of_perturbing}
	\end{proposition}
	
	\begin{proof}
	    Since $\tilde{U}$ is a GLM, we have $\nabla E(\tilde{U}) =0$. For any $\delta U \in \mathcal{W}^k(U)$, we have
		$\langle \delta U, H(\tilde{ U})\delta U \rangle = 0$.
		We employ the Taylor expansion of $E(U)$ at $\tilde{U}$,
	    \begin{align}
		    E(\tilde{U} + \delta U) & = E(\tilde{U}) + \nabla E(\tilde{U})\delta U+ \frac{\langle \delta U, H(\tilde{ U})\delta U \rangle}{2!} + o(\|\delta U\|_{\ell^2}^2) \label{eq:translational_invariance}\\
		& = E(\tilde{U}) + o(\|\delta U\|_{\ell^2}^2),\notag
	    \end{align}
	    where $\delta U \in \mathcal{W}^k(U)$ satisfies $\|\delta U\|_{\ell^2} < \epsilon $. Then we have
		$|E(\tilde{U} + \delta U) - E(\tilde{U})| \leq  o(\epsilon^2)$. The proof is completed.
	\end{proof}
	
	\Cref{prop:stability_of_perturbing} indicates that the energy remains almost unchanged by a minor perturbation of the nullspace for the GLM. For phase transitions with translational invariance, existing surface-walking methods that search for saddle points along the eigenvector $\bm{v}$ linked to the smallest eigenvalue, such as the gentlest ascent dynamics and dimer-type methods, could be challenging to escape from the attraction basin\,\cite{yin2021transition}. Since $\bm{v} \in \mathcal{W}^k(\tilde{U})$, a slight movement along $\bm{v}$ does not result in an increase in energy value at the GLM $\tilde{U}$.
	In order to avoid being trapped in a basin, we introduce the NPSS method in \Cref{sec:method}, which ensures that the search direction is always upward. 
	%We present the NPSS method in \Cref{sec:method}, and analyze how to deal with such degenerate problems.

    %Due to the stability of energy in $\mathcal{W}^k(\bar{U})$, if the vector in $\mathcal{W}^k(\bar{U})$ is taken as ascent direction, it will not escape from the basin.	
%	\newpage

	\section{Nullspace-preserving saddle search method}
	\label{sec:method}
    In this section, we propose the NPSS method to find the index-1 GSP starting from a GLM on the energy surface. The NPSS method consists of two stages. In Stage I, for any state $U$ in the GQR of GLM, the NPSS method ascends in segments based on the variation of the nullspace. The ascent direction is maintained orthogonal to the nullspace of the initial state in each segment. The initial state in each segment is updated when the nullspace of the current state differs significantly from that of the initial state. These operations ensure the effectiveness of the ascent direction and reduce the costs of updating the nullspace at each step, enabling a quick escape from the basin. Once a negative eigenvalue of the Hessian appears, indicating that the system has escaped from the GQR. In Stage II, the NPSS method ascends along the ascent subspace $\mathcal{V}$ spanned by the eigenvector corresponding to the negative eigenvalue and descends along the orthogonal complement subspace $\mathcal{V}^{\bot}$. Then the system can gradually converge to the index-1 GSP. The mathematical description of the NPSS method is provided below.
    
	\subsection{Stage I : Escape from basin}
	We consider a state $U$ in the GQR of $U^{(0)}$, whose Hessian's subspaces contain $\mathcal{W}^k(U)$ and $\mathcal{W}^s(U)$. 
	An easier way to escape the basin is to follow the gentlest ascent direction, i.e., the eigenvector $\bm{v}$ corresponding to the smallest positive eigenvalue. We define an ascent space $\mathcal{V} = {\rm span}\{\bm{v}\}$ and denote its orthogonal complement as $\mathcal{V}^{\bot}$. The negative gradient is $T(U) = -\nabla E(U)$, where $\nabla E(U)$ represents the steepest ascent direction at $U$. Let $\mathcal{P}_{\mathcal{V}}$ be the orthogonal projection operator onto the subspace $\mathcal{V}$. To efficiently escape from the basin, we take $-\mathcal{P}_{\mathcal{V}} T(U)$ as the ascent direction in $\mathcal{V}$ and $\mathcal{P}_{\mathcal{V}^{\bot}} T(U)$ as the descent direction in the orthogonal complement $\mathcal{V}^{\bot}$. Escaping from the basin can be achieved by the gradient-type dynamics
	\begin{align}
		\dot{U}&=\beta_1\left(-\mathcal{P}_{\mathcal{V}} T(U)\right) + \beta_2 \mathcal{P}_{\mathcal{V}^{\bot}} T(U) \\
		&=\beta_1\left(-\mathcal{P}_{\mathcal{V}} T(U)\right) + \beta_2\left(I -\mathcal{P}_{\mathcal{V}} \right)T(U),\notag
	\end{align}
	where $\beta_1$ and $\beta_2$ represent positive relaxation constants. For simplicity, we set $\beta_1 = \beta_2 = \beta$, and can obtain a modified direction, resulting in a modified gradient-type dynamic,
	\begin{align}
		\beta^{-1} \dot{U}=(I-2 \mathcal{P}_{\mathcal{V}} )T(U).
		\label{eq:update-u}
	\end{align}
	
		\begin{figure}[!hbpt]
		\centering
		\includegraphics[scale=0.08]{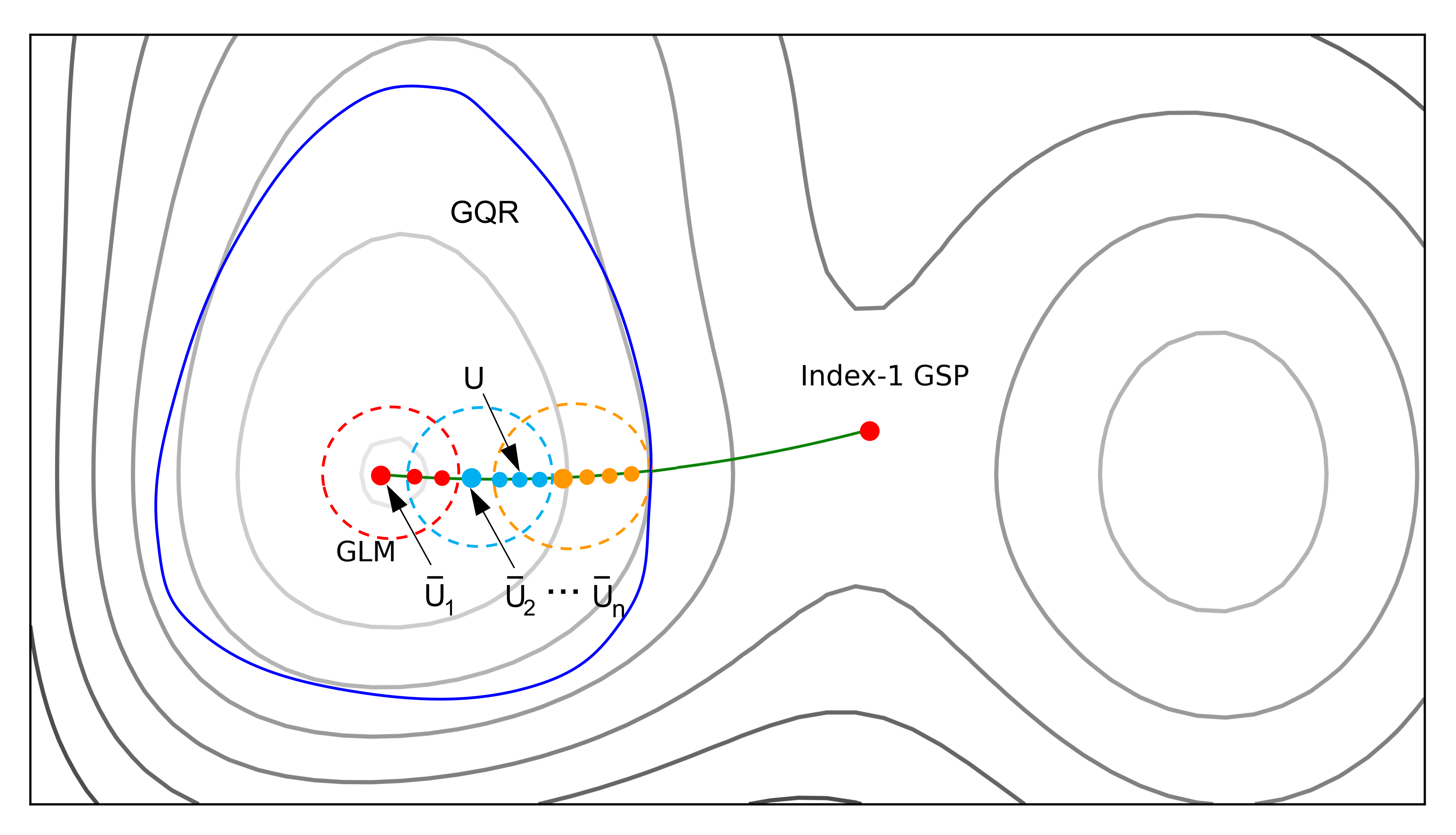}
		\caption{Illustration of an energy surface. The evolution path is shown as a green curve. The GLM and index-1 GSP are represented by red points on the evolution path. The boundary of GQR is marked by the blue curve. $\{\bar{U}_i\}_{i=1}^n$ are the initial states of each segment in the evolution path. The shallow circles at $\bar{U}_i$ represent a neighborhood satisfying $\|\sin \Theta(\mathcal{W}^k(U), \mathcal{W}^k( \bar{U}_i))\|_2 < 1$ for any state $U$ within the GQR.}
		\label{fig:PES}
	\end{figure}
	
	The key to escaping from the basin is to determine the ascent subspace $\mathcal{V}$. For this degenerate problem, directly minimizing the Rayleigh quotient $\langle \bm{v}, H(U) \bm{v} \rangle / \langle \bm{v}, \bm{v} \rangle $ can only obtain eigenvectors corresponding to zero eigenvalues. If the vector in the nullspace is mistakenly taken as the ascent direction, the system will fail to escape the basin. To avoid the effects of nullspace, we should theoretically impose the constraint $\bm{v}\bot \mathcal{W}^k(U)$.
    However, accurately calculating $\mathcal{W}^k(U)$ for varying $U$ is expensive in numerical experiments. An intuitive idea is to replace $\mathcal{W}^k(U)$ directly with $\mathcal{W}^k(U^{(0)})$, but the error may not be effectively controlled as the system state evolves. 
    
    A better way is to consider nullspace replacement in segments. As shown in \Cref{fig:PES}, $\{\bar{U}_i\}_{i=0}^n$ denote the initial states of $n$ segments. For the $i$-th segment with the initial state $\bar{U}_i$, we use its nullspace $\mathcal{W}^k(\bar{U}_i)$ to replace $\mathcal{W}^k(U)$ for the subsequent state $U$ within this segment. The difference between $\mathcal{W}^k(U)$ and $\mathcal{W}^k(\bar{U}_i)$ is measured by the principal angles. When $\|\sin \Theta (\mathcal{W}^k(U), \mathcal{W}(\bar{U}_i))\|_2 = 1$, we use the current state $U$ to update $\bar{U}_{i+1}$ and start a new segment. Without loss of generality, we consider $\bar{U} \in \{\bar{U}_i\}_{i=0}^n$ as the initial state of any segment.  Then we can solve the optimization problem \cref{eq:min-v} to obtain $\bm{v}\in \mathcal{V}$ in this segment.
	\begin{align} \label{eq:min-v}
		\min _{\bm{v}}\langle \bm{v}, H(U) \bm{v} \rangle 
		\quad \text {s.t.} \quad \left\{ \begin{array}{l}
			\langle \bm{v}, \bm{v}\rangle = 1,\\
			\bm{v} \bot \mathcal{W}^k(\bar{U}).
		\end{array} \right.
	\end{align} 
	We propose a lemma to show that the optimization problem \cref{eq:min-v} can ensure effectively obtaining an ascent direction.
	
	\begin{lemma}\label{le:find_direction}
    The solution $\bm{v}$ of \cref{eq:min-v} exists an ascent component for any state $U$ in the neighborhood of $\bar{U}$, satisfying 
    \begin{align}
    \|\sin \Theta(\mathcal{W}^k( U), \mathcal{W}^k( \bar{U} )\|_2 < 1.
    \end{align}
	\end{lemma}
	
	\begin{proof}
		 For a state $U$ in the neighborhood of $\bar{U}$, there is
		\begin{align}
		\bbR^M = \mathcal{W}^k(U) \oplus \mathcal{W}^s(U), \label{eq:space_R}
		\end{align}
		where $\mathcal{W}^k(U)$ and $\mathcal{W}^s(U)$ can be spanned by unit orthogonal eigenvectors $\{\bm{v}^k_i\}_{i=1}^{l_k}$ and $\{\bm{v}^s_i\}_{i=1}^{l_s}$, respectively.
        
        Let $P_{\mathcal{W}^k(\bar{U})}$ and $P_{\mathcal{W}^k(U)}$
        be the orthogonal projection operators of $\mathcal{W}^k(\bar{U})$ and $\mathcal{W}^k(U)$, respectively.
        Based on the relationship between the projection operators and $\sin \Theta$ as defined in \cref{eq:sin_theta}, we have
        \begin{align}
            \|P_{\mathcal{W}^k(U)} (I - P_{\mathcal{W}^k(\bar{U})} )\|_2
            &= \|\sin \Theta(\mathcal{W}^k( U), \mathcal{W}^k( \bar{U}))\|_2.\notag
        \end{align}
        For any $\bm{v} \in \bbR^M$ satisfying $\langle \bm{v}, \bm{v}\rangle = 1$ and $\bm{v} \bot \mathcal{W}^k(\bar{U})$, 
        we have 
        \begin{align}
        \bm{v} = (I - P_{\mathcal{W}^k(\bar{U})})\bm{v}.
        \end{align}
        Then\,\cref{eq:project_v} can be obtained,
        \begin{align}\label{eq:project_v}
            \|P_{\mathcal{W}^k(U)} \bm{v}\|_{2} 
            &=\|P_{\mathcal{W}^k(U)} (I - P_{\mathcal{W}^k(\bar{U})} ) \bm{v}\|_{2}\\ \notag
            &\leq \|P_{\mathcal{W}^k(U)} (I - P_{\mathcal{W}^k(\bar{U})} )\|_2\\ \notag
            &= \|\sin \Theta(\mathcal{W}^k( U), \mathcal{W}^k( \bar{U}))\|_2 \leq 1 . \notag
        \end{align}
        Denote $C_{\Theta} = \|\sin \Theta(\mathcal{W}^k(U), \mathcal{W}^k( \bar{U}))\|_2$ is the sine of the largest principal angle between $\mathcal{W}^k(U)$ and $\mathcal{W}^k(\bar{U})$. $C_{\Theta} = 1$ holds only when the largest principal angle is $\pi/2$.

        According to \cref{eq:space_R}, $ P_{\mathcal{W}^k(U)} \bm{v}$ can be written as 
        \begin{align}
            P_{\mathcal{W}^k(U)} \bm{v} &= c_{l_s + 1} \bm{v}_{1}^{k} +\cdots+ c_{M} \bm{v}_{l_k}^{k} ,\notag
        \end{align}
        where $(c_{l_s + 1})^2+ \cdots +(c_{M})^2 \leq C_{\Theta}^2$.
        Similarly, $P_{\mathcal{W}^s(U)}\bm{v}$ can be written as  $P_{\mathcal{W}^s(U)}\bm{v} = c_1 \bm{v}_{1}^{s} + \cdots + c_{l_s} \bm{v}_{l_s}^{s}$, then
        \begin{align}
            \bm{v} &= P_{\mathcal{W}^s(U)}\bm{v} + P_{\mathcal{W}^k(U)} \bm{v}\\
            &= c_1\bm{v}_{1}^{s} + \cdots + c_{l_s} \bm{v}_{l_s}^{s} + c_{l_s + 1} \bm{v}_{1}^{k} + \cdots + c_{M} \bm{v}_{l_k}^{k},\notag
        \end{align}
        with $(c_{1})^2 + \cdots + (c_{l_s})^2 + (c_{l_s + 1})^2+ \cdots +(c_{M})^2 = 1$. Then we have $(c_{1})^2 + \cdots + (c_{l_s})^2 \geq 1-C_{\Theta}^2$. Let $\lambda^s_i$ be eigenvalues corresponding to basis vectors $\bm{v}^s_i$, $i=1, 2,\cdots,l_s$, then
        \begin{align}
            \langle \bm{v}, H(U) \bm{v}\rangle &= (c_{1})^2\lambda^s_{1} +(c_{2})^2\lambda^s_{2} + \cdots + (c_{l_s})^2\lambda^s_{l_s} + (c_{l_s + 1})^2 0 + \cdots + (c_{M})^2 0\\\notag
            &= (c_{1})^2\lambda^s_{1} +(c_{2})^2\lambda^s_{2} + \cdots + (c_{l_s})^2\lambda^s_{l_s}.\notag
        \end{align}
        Without loss of generality, let $\lambda_1^{s}$  be the minimum positive eigenvalue corresponding to $\bm{v}_1^s$, then
        \begin{align}
        \langle \bm{v}, H(U) \bm{v}\rangle &\geq ((c_1)^2 + \cdots + (c_{l_s})^2)\lambda_1^{s}\geq (1-C_{\Theta}^{2})\lambda_1^{s}. 
        \end{align} 
		Therefore, the minimum of optimization problem \cref{eq:min-v} is 
		\begin{align} \label{eq:lambda_min}
		 \lambda_{\min} = (1-C_{\Theta}^{2})\lambda_1^{s},
		\end{align}
		and the corresponding eigenvector $\bm{v}$ is
        \begin{align}
            \bm{v}=\sqrt{1-C_{\Theta}^{2}} \bm{v}_{1}^{s}+C_{\Theta} \bm{v}^{k} ,\label{eq:ascent_direction}
        \end{align}
        where $\bm{v}^{k} \in \mathcal{W}^k(U)$ satisfies $\|\bm{v}^{k}\|_2 = 1$. When $C_{\Theta} < 1$, $\bm{v}$ always can provide an ascent component $\sqrt{1-C_{\Theta}^{2}} \bm{v}_{1}^{s}$ to push the system escape from the basin. The proof is completed.
        \end{proof}

        In the following, we discuss how $C_{\Theta}$ determines the acquisition of the ascent direction.
        
        (i). When $ 0\leq C_{\Theta} < 1$, the component $\sqrt{1-C_{\Theta}^{2}} \bm{v}_{1}^{s}$ can provide an efficient ascent direction. We can obtain the ascent direction by solving the optimization problem \cref{eq:min-v} with gradient-type dynamics,
        \begin{align}
		\left\{\begin{array}{l}
			\zeta^{-1}\dot{\bm{v}}=-H(U) \bm{v}+\left\langle \bm{v}, H(U) \bm{v}\right\rangle \bm{v} + \sum_{i=1}^{l_k}\langle \bar{\bm{v}}^{k}_{i}, H(U) \bm{v} \rangle  \bar{\bm{v}}^{k}_{i} ,\\
			\bm{v} = \bm{v}/\|\bm{v}\|_{\ell^2}.
		\end{array}\right. 
		\label{eq:update-v}
	\end{align}
	where $\zeta>0$ is a relaxation parameter, $\mathcal{W}^k( \bar{U}) = {\rm span} \{ \bar{\bm{v}}_{1}^{k}, \bar{\bm{v}}_{2}^{k},\cdots, \bar{\bm{v}}_{l_k}^{k}\}$ and $\{\bar{\bm{v}}^{k}_{i}\}_{i=1}^{l_k}$ are unit orthogonal bases. 
    
    (ii). When $C_{\Theta} = 1$, it means that the largest principal angle between  $\mathcal{W}^k(U)$ and $\mathcal{W}^k( \bar{U}))$ is $\pi/2$. There is a base in $\mathcal{W}^k(U)$ that is orthogonal to $\mathcal{W}^k( \bar{U})$ and $\bm{v} = \bm{v}^{k}$, which no longer provides an ascent direction. Therefore, we update $\bar{U} = U$ and $\mathcal{W}^k(\bar{U}) = \mathcal{W}^k(U)$, and start a new segment to climb upwards by solving\,\cref{eq:update-v}. These operations can ensure that an efficient ascent direction always be obtained.  
    
    \begin{remark}
     In numerical experiments, computing directly $C_{\Theta}$ is expensive. Alternatively, we adopt $\lambda_{\min}$ in \cref{eq:lambda_min} as an efficient indicator to determine the update of nullspaces. For any $U$ within the GQR of the GLM, we have $\lambda_1^{s} > 0$, hence conditions $\lambda_{\min} \to 0$ and $C_{\Theta} \to 1$ are equivalent. 
     Therefore, we consider updating $\bar{U}$ and the nullspace only if
    \begin{align}
        \lambda_{\min} < \epsilon_{\lambda},
    \end{align}
    where $\epsilon_{\lambda} > 0$ is a threshold. 
\end{remark}

	\subsection{Stage II : Search for index-1 GSP}
	Once the minimum eigenvalue of $H(U)$ becomes negative, it means that $U$ escapes the GQR of the GLM. Since the index-1 GSP has only one negative eigenvalue, it is the maximum in the direction of the eigenvector corresponding to the negative eigenvalue, and the minimum in the remaining directions.  Thus, to search the index-1 GSP, we span the ascent space $\mathcal{V}$ by the eigenvector $\bm{v}$ corresponding to the minimum eigenvalue, i.e., the negative eigenvalue. Then, we search for a local maximum on the linear manifold $U+ \mathcal{V} =\{U +\bm{v}| \forall \bm{v} \in \mathcal{V}\}$, and a local minimum on $U+ \mathcal{V}^{\bot}$\,\cite{agusti2011affine}.
	This process can be transformed into a minimax optimization problem
	\begin{align}
		 \min _{U^d \in \mathcal{V}^{\bot}}  \max_{U^a \in \mathcal{V}} ~ E\left(U^d + U^a \right), \label{eq:minimax}
	\end{align}
	where $U^a$ and $U^d$ are the projections of $U$ on $\mathcal{V}$ and $\mathcal{V}^{\bot}$. Many methods can solve this problem, such as the Newton method, and gradient-type dynamics.
	The latter can be written as
	\begin{align}
		\dot{U} &=\beta_1 (-\mathcal{P}_{\mathcal{V}}T(U)) +\beta_2\mathcal{P}_{\mathcal{V}^{\bot}}T(U) \label{eq:update-u2}\\
		& =\beta_1 (-\mathcal{P}_{\mathcal{V}}T(U)) +\beta_2(I-\mathcal{P}_{\mathcal{V}})T(U)\notag\\
		& = \beta(I - 2\mathcal{P}_{\mathcal{V}})T(U),\notag
	\end{align}
	where $\beta_1 = \beta_2 = \beta$ are positive relaxation constants, $-\mathcal{P}_{\mathcal{V}}T(U)$ is the ascent direction along $\mathcal{V}$ and $\mathcal{P}_{\mathcal{V}^{\bot}}T(U)$ is the descent direction along $\mathcal{V}^{\bot}$. The dynamics \cref{eq:update-u2} have been rigorously proven that its stable fixed points correspond to index-1 GSPs\,\cite{e2011gentlest}. Furthermore, the subspace $\mathcal{V} = {\rm span}\{\bm{v}\}$ can be obtained by directly minimizing the Rayleigh quotient,
	\begin{align}\label{min-v2}
		\min _{\bm{v}}\langle \bm{v}, H(U) \bm{v} \rangle  
	    \quad \text {s.t.} \quad	\langle \bm{v}, \bm{v} \rangle = 1.
	\end{align} 
	The dynamic form is 
	\begin{align}
		\left\{\begin{array}{l}
			\zeta^{-1}\dot{\bm{v}}=-H(U) \bm{v} +\left\langle \bm{v}, H(U) \bm{v} \right\rangle \bm{v} ,\\
			\bm{v} = \bm{v}/\|\bm{v} \|_{\ell^2}.
		\end{array}\right. 
		\label{eq:update-v2}
	\end{align}
	
	\Cref{alg:NPSS} summarizes the NPSS method.
	\begin{algorithm}
		\caption{NPSS method}
		\label{alg:NPSS}
		\begin{algorithmic}
			\STATE Stage I : Escape from the basin
			\STATE Input : $U^{(0)}$, $\bm{v}^{(0)} \in \mathbb{R}^M$ satisfying $\bm{v}^{(0)} \perp \mathcal{W}^k(U^{(0)})$,  $\mathcal{W}^k(U^{(0)}) \subset \mathbb{R}^M$, $\epsilon_{\lambda},\epsilon_{T}$
			\STATE Set $n =1$, $U^{(n)} = U^{(0)}+\xi \bm{v}^{(0)}$, compute $T(U^{(n)})$
		    \STATE Set $\bar{U} = U^{(0)}$, $\mathcal{W}^k(\bar{U}) =\mathcal{W}^k(U^{(0)})$
			\REPEAT 
			\STATE Update $U^{(n+1)}$ as $\cref{eq:update-u}$
			\STATE Calculate $\bm{v}^{(n+1)}$ as $\cref{eq:update-v}$
			\STATE Calculate $T(U^{(n+1)})$
			\STATE $n:= n+1$
			\STATE If $\langle \bm{v}, H(U^{(n)}) \bm{v} \rangle <\epsilon_{\lambda}$
			\STATE  \qquad Update  $\bar{U} = U^{(n)}$, $\mathcal{W}^k(\bar{U}) =\mathcal{W}^k(U^{(n)})$ 
			\UNTIL{$ \langle  \bm{v}^{n},
				H(U^{(n)})\bm{v}^{n}\rangle < 0$}
			\STATE Output :  $U^{(n)}$, the negative gradient $T(U^{(n)})$ and the ascent direction $\bm{v}^{(n)}$
			\vspace{0.1cm} \hrule \vspace{0.1cm}
			\STATE Stage II : Search for index-1 GSP
			\REPEAT 
			\STATE Update $U^{(n+1)}$ as $\cref{eq:update-u2}$
			\STATE Calculate $\bm{v}^{(n+1)}$ as $\cref{eq:update-v2}$
			\STATE Calculate $T(U^{(n+1)})$
			\STATE $n:= n+1$;
			\UNTIL{$\|T(U^{(n)})\|_{\ell^2} < \epsilon_T$}
			\STATE Output: the index-1 GSP $U^{(n)}$
		\end{algorithmic}
	\end{algorithm}
	
	Consider the LP model as an example. In \Cref{alg:NPSS}, the gradient dynamics \cref{eq:update-u} and \cref{eq:update-u2} for updating $U^{(n+1)}$ have the same form in both stages. We discretize them with the semi-implicit scheme as \cref{discretedynamic_u}.
	\begin{align}
	    \beta^{-1}(U^{(n+1)}-U^{(n)}) =& -\mathcal{F}^{-1}_M D \mathcal{F}_M U^{(n+1)} +\varepsilon  U^{(n+1)} \label{discretedynamic_u}\\
	    & + \alpha  (U^{(n)})^2- (U^{(n)})^3 - 2 \langle \bm{v}^{(n)},T(U^{(n)})\rangle \bm{v}^{(n)},\notag
	\end{align}
    where $ T(U^{(n)}) = -\mathcal{F}^{-1}_M D \mathcal{F}_M U^{(n)} +\varepsilon  U^{(n)} + \alpha (U^{(n)})^2- (U^{(n)})^3$.  For the the ascent direction $\bm{v}^{(n+1)}$, we adopt the forward Euler method to discrete $\cref{eq:update-v}$ as \cref{discretedynamic_v1} and $\cref{eq:update-v2}$ as \cref{discretedynamic_v2},
    \begin{align}
	\zeta^{-1}(\bm{v}^{(n+1)} - \bm{v}^{(n)})=&-H(U^{(n+1)}) \bm{v}^{(n)} +\langle \bm{v}, H(U^{(n+1)}) \bm{v}^{(n)}\rangle \bm{v}^{(n)} \label{discretedynamic_v1}\\
	&+ \sum_{i=1}^{l_k}\langle \bar{\bm{v}}^{k}_{i}, H(U^{(n+1)}) \bm{v}^{(n)} \rangle  \bar{\bm{v}}^{k}_{i},  \notag \\
	 \zeta^{-1}(\bm{v}^{(n+1)} - \bm{v}^{(n)})=& -H(U^{(n+1)}) \bm{v}^{(n)} +\langle \bm{v}^{(n)}, H(U^{(n+1)}) \bm{v}^{(n)} \rangle \bm{v}^{(n)}.\label{discretedynamic_v2}
	\end{align}

    \subsection{An illustrative example} 
\begin{figure}[H]
		\centering
		\includegraphics[width=0.95\textwidth]{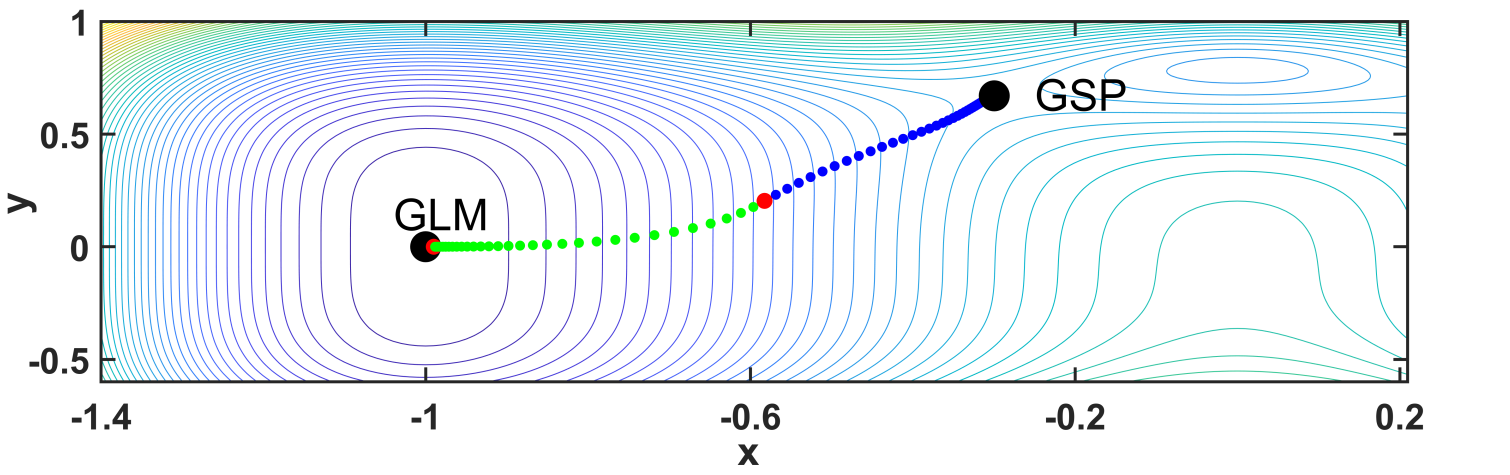}
		\caption{Evolution path from the GLM to index-1 GSP using the NPSS method. Black dots represent the GLM and index-1 GSP, respectively. Green and blue dots represent the evolution path before and after reaching the boundary of GQR, and red dots indicate the points where the nullspace is recalculated.}
		\label{fig:NPSS_functioncase}
	\end{figure}

    To better understand the key idea of the NPSS method, we consider a function case
    $f(x, y) = x^4 + y^4 - 2x^2 - e^{-6(x^2 + (y-1)^2)}$, which has a degenerate GLM at $(-1, 0.0249)$ and an index-1 GSP at $(-0.2996, 0.6698)$, as shown in \Cref{fig:NPSS_functioncase}.
    
    We use the NPSS method to search for the index-1 GSP from the GLM. We set $\beta = 0.01$ and $\zeta = 0.01$ as the iteration steps, $\|\nabla f\|_{\ell^2}<$1e-7 as the convergence condition and  $\epsilon_{\lambda} = 0.05$ as the threshold for updating the nullspace.  \Cref{fig:NPSS_functioncase} demonstrates that the NPSS method can efficiently push the system out of the GLM basin and converge to the index-1 GSP. Notably, the NPSS method requires only a few times of nullspace updates to escape the basin, saving computational time by avoiding nullspace calculations at each step.

	\section{Numerical examples}\label{sec:result}
	
	In this section, to study phase transitions with translational invariance, we apply the gradient descent (GD) method to find the GLM and the NPSS method to search for index-1 GSP\,\cite{yin2021transition,zhang2008efficient,jiang2020efficient}. Several numerical experiments based on LB and LP models are presented: (i). Nucleation from DIS to ordered phases (see \Cref{nucleation_from_dis}); (ii). Phase transitions between crystals (see \Cref{crystal_and_crystal}); (iii). Phase transitions between quasicrystal and crystal (see \Cref{quasicrystal_and_crystal}). Furthermore, we demonstrate the efficiency and reliability of the NPSS method in comparison with the HiSD method in \Cref{quasicrystal_and_crystal}. 
	
	For these experiments, we use $[0, 2\pi L)^2$ in 2D and $[0, 2\pi L)^3$ in 3D as the computational domains, and $N$ spatial discretization points in each dimension. The convergence condition for the saddle point search is set to $\|F(U)\|_{\ell_2}<$ 1e-7. All methods are implemented in C++ programming language and run on the same workstation, 12th Gen Intel (R) Core (TM) i7-12700KF, 16 GB memory, under Linux.

	\newpage
	\subsection{Nucleation from disordered phases}\label{nucleation_from_dis}

\begin{figure}[!hbpt]
    \begin{minipage}[t]{1\linewidth}
    \centering
    \includegraphics[width=5in,height=2.3in]{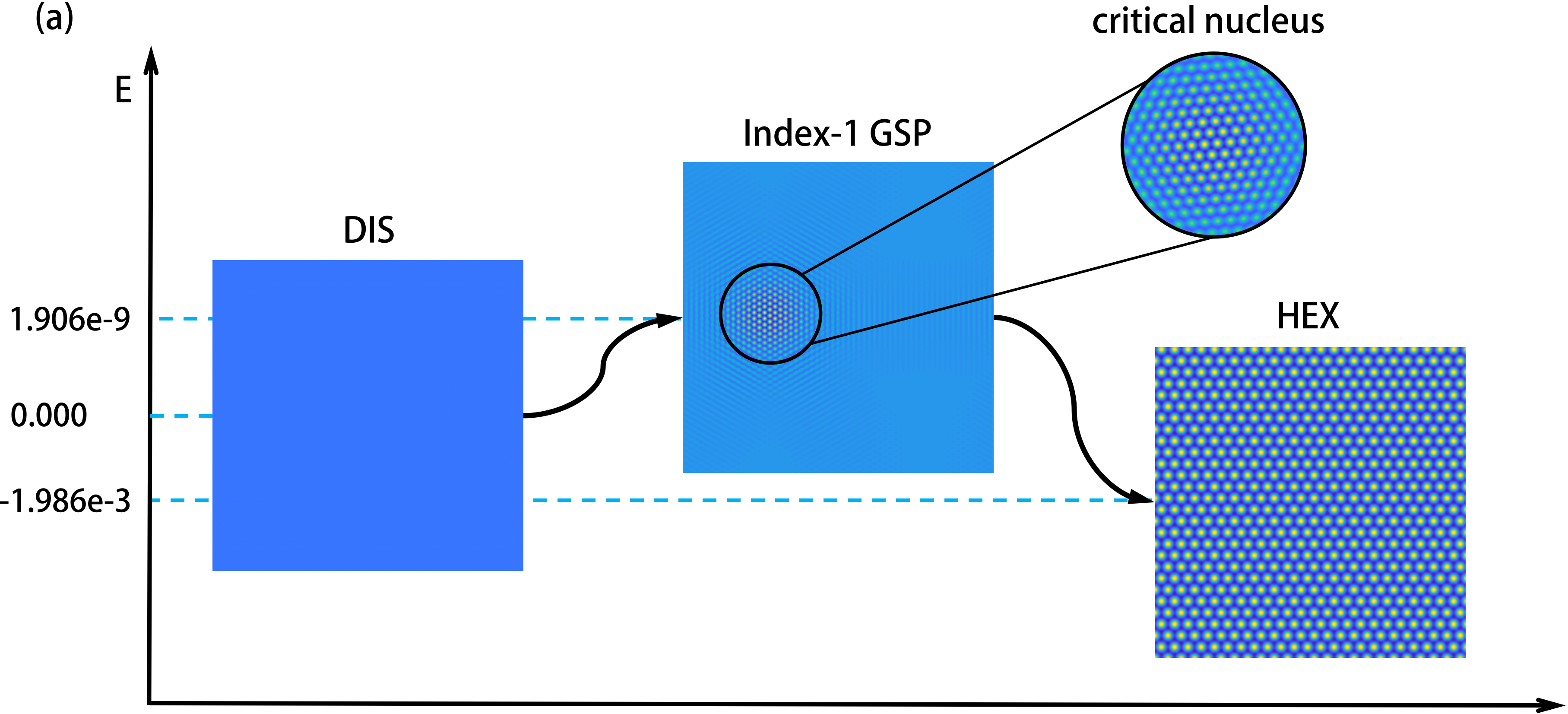}\\
    \end{minipage}%
    \\
    \begin{minipage}[t]{1\linewidth}
    \centering
    \includegraphics[width=5.07in,height=2.35in]{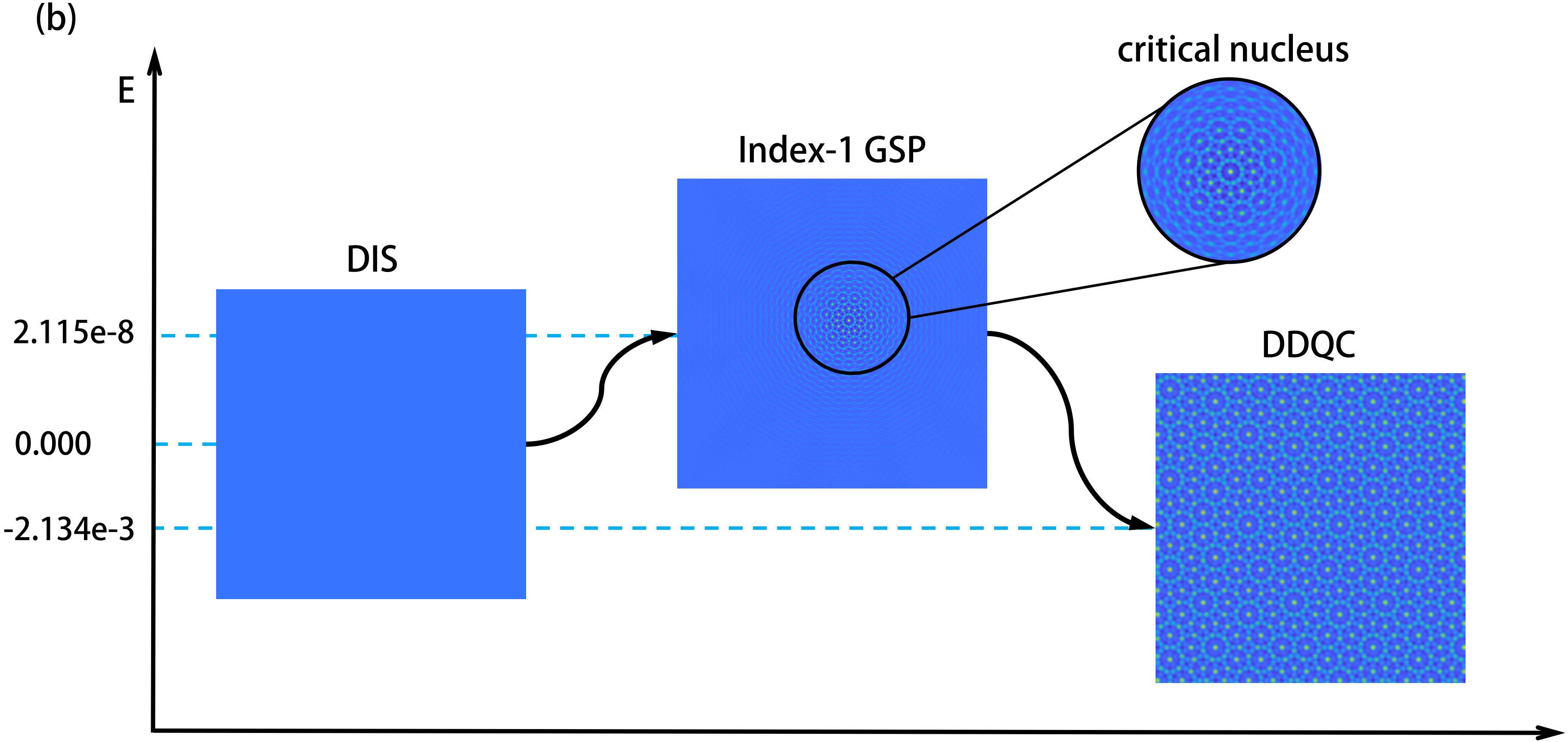}\\
    \end{minipage}%
    \caption{(a). Transition path from DIS to HEX computed by the NPSS method in LB model with $\tau = 0.001$, $\gamma = 0.5$, where $L = 60$, $N = 900$. (b). Transition path from DIS to DDQC computed by the NPSS method in LP model with $\varepsilon = -0.01$, $\alpha = 0.95$, where $L = 112$, $N = 1024$. The $x$-axis represents the evolution direction of phase transitions.}
    \label{fig:DIStoOrdered}
    \end{figure}

    In this subsection, we use the NPSS method to study the nucleations of the crystal and quasicrystal from DIS. At parameters in both two cases (see \Cref{fig:DIStoOrdered}), the Hessian of DIS has only positive eigenvalues, thus the NPSS method only needs to climb along the eigenvector direction corresponding to a positive eigenvalue. We show phase transition paths from DIS to HEX using the LB model and to DDQC within the LP model, demonstrating the power of the NPSS method.

	\subsection{Phase transition between crystals} \label{crystal_and_crystal}
	
\begin{figure}[H]
		\centering
		\includegraphics[width=5.07in,height=2.35in]{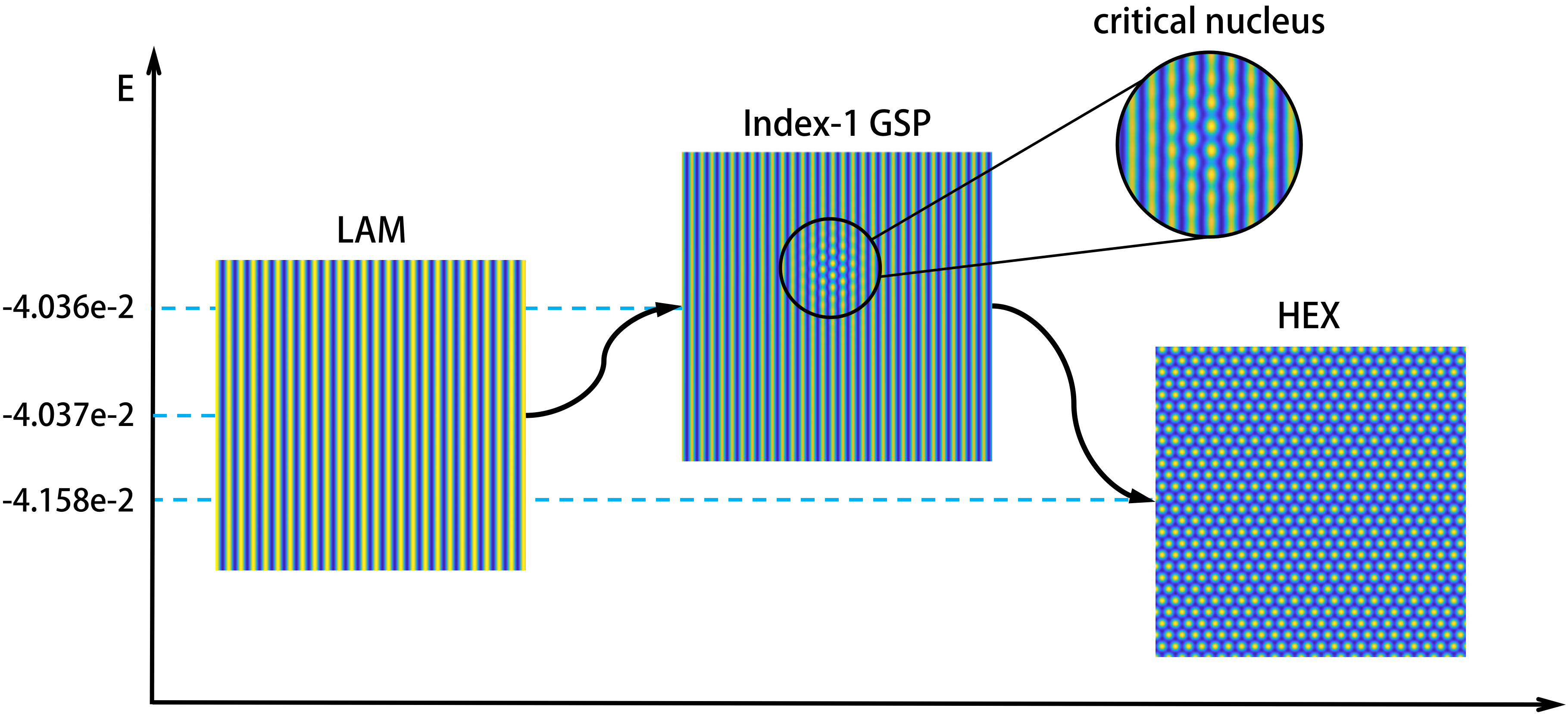}
		\caption{Transition path from LAM to HEX computed by the NPSS method in the LB model with $\tau = -0.2$, $\gamma = 0.28$, where $L = 60 $, $N = 900$.}
		\label{fig:LAMtoC6}
	\end{figure}

	For the LB model, we compute the phase transition between LAM and HEX using the NPSS method. At $\tau = -0.2$ and $\gamma = 0.28$, LAM is a metastable state with $E$ = -4.037e-2, while HEX is a relatively stable state with $E$ = -4.158e-2. As shown in \Cref{fig:LAMtoC6}, a circular critical nucleus of HEX with $E$ = -4.036e-2 is found in the transition path from LAM to HEX. A similar phenomenon has also been obtained by the string method \cite{li2013nucleation}. 
	
\begin{figure}[!ht]
		\centering
		\includegraphics[width=1\textwidth]{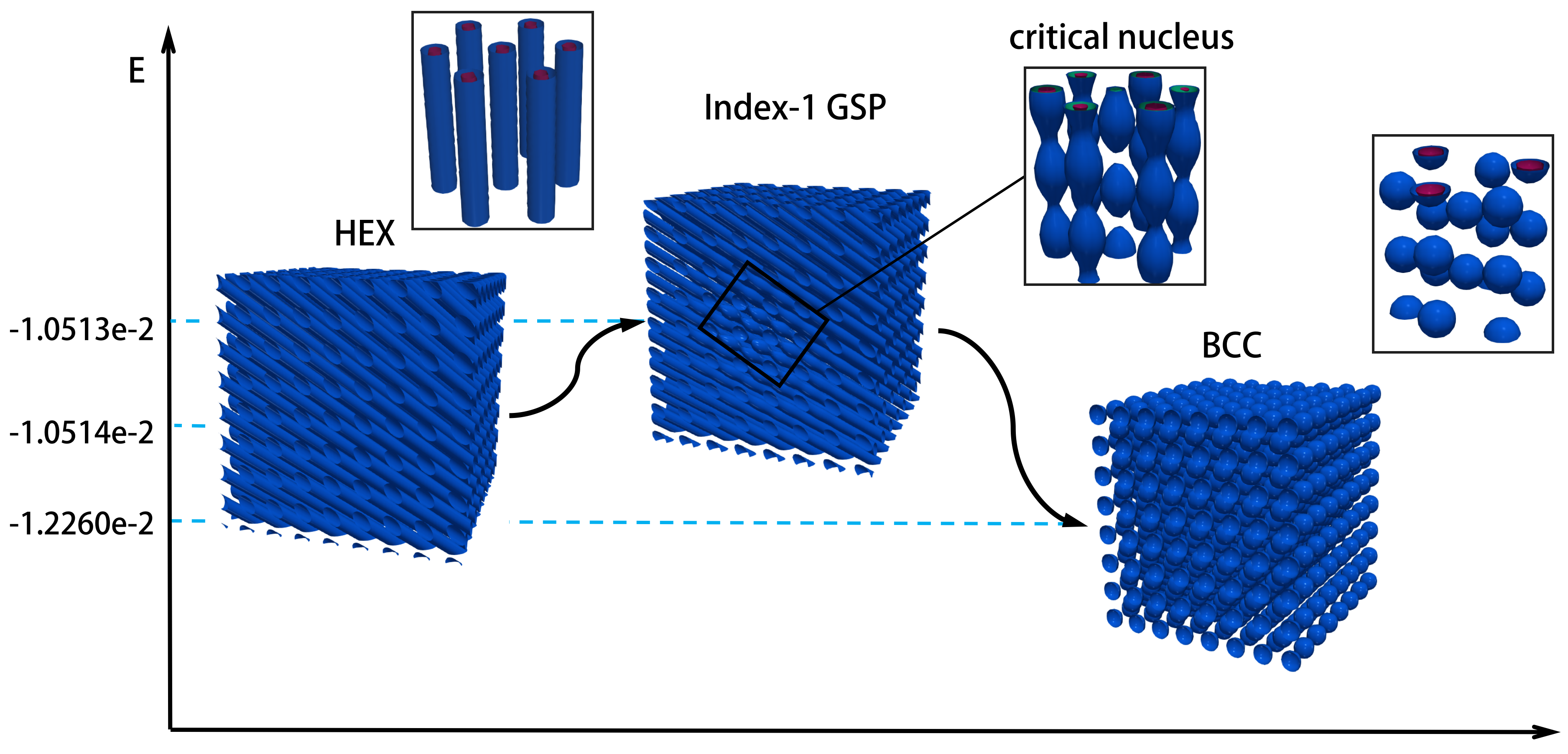}
		\caption{Phase transition from 3-dimensional HEX to BCC computed by the NPSS method in the LB model with $\tau = -0.008$, $\gamma = 0.7$, where $L = 20\sqrt{2}$ and $N = 224$. The three insert images show the local structures along the body diagonal line.}
		\label{fig:HEXtoBCC}
	\end{figure}

	We also use the NPSS method to study 3D ordered phase transitions and locate the transition state in the LB model. At $\tau = -0.008$ and $\gamma = 0.7$, we study the phase transition from HEX with $E$ = -1.0514e-2 to BCC with $E$ = -1.2260e-2. The complete transition path is shown in \Cref{fig:HEXtoBCC}. From the insert images, the columnar structures in HEX partially deform to form a critical nucleus with $E$ = -1.0513e-2 and eventually evolve into the spherical BCC structure. This result is consistent with the experimental observation in morphological evolution\,\cite{sota2013directed}.
	
	\subsection{Phase transition between quasicrystals and crystals}\label{quasicrystal_and_crystal} 
	We apply the NPSS method to study the phase transition between DDQC and HEX based on the LP model and compare the computational efficiency with the HiSD method. Let $L= 112$ and $N=1024$. We choose $U^* = U^{(0)}+\xi \bm{v}$ as the initial position to search for saddle points, where $U^{(0)}$ is a metastable state that can be obtained by the gradient flow with a proper initial value\,\cite{jiang2015stability}. $\bm{v}$ is the perturbation direction and is usually set as the normalized eigenvector of the smallest positive eigenvalue, which can be computed by the locally optimal
    block preconditioned conjugate gradient (LOBPCG) method\,\cite{knyazev2001toward}. $\xi$ is a positive constant used to push the system away from the minimum. \Cref{tab:parameters} gives the specific values used in this experiment.
    
	\begin{table}[!hbpt]
	    \centering
	    \footnotesize{
	    \caption{The perturbation $\xi$ at the initial position during the saddle point search for different phase transitions.}
	    \begin{tabular}{|c|c|c|c|c|}
	    \hline
	    Stage & DDQC $\rightarrow$ LQ & LQ $\rightarrow$ HEX & HEX $\rightarrow$ LQ & LQ $\rightarrow$ DDQC \\
	       \hline
	       $\xi$ & 0.1 & 0.15& 0.01& 0.1\\
	       \hline
	    \end{tabular}
	    }
	    \label{tab:parameters}
	\end{table}
	
\begin{figure}[!hbpt]
        \begin{minipage}[t]{1\linewidth}
        \centering
        \includegraphics[width=5in,height=2.1in]{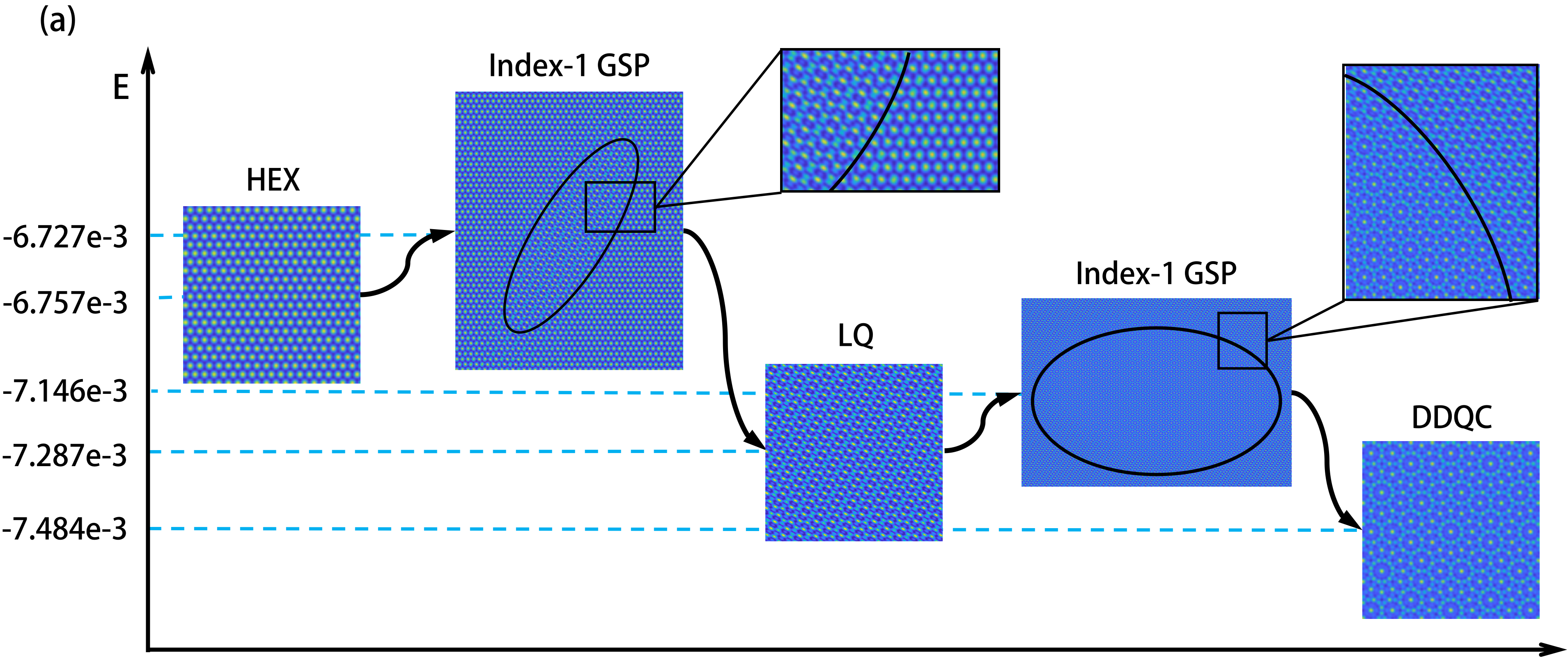}\\
        \end{minipage}%
        \\
       \begin{minipage}[t]{1\linewidth}
        \centering
        \includegraphics[width=5.05in,height=2.1in]{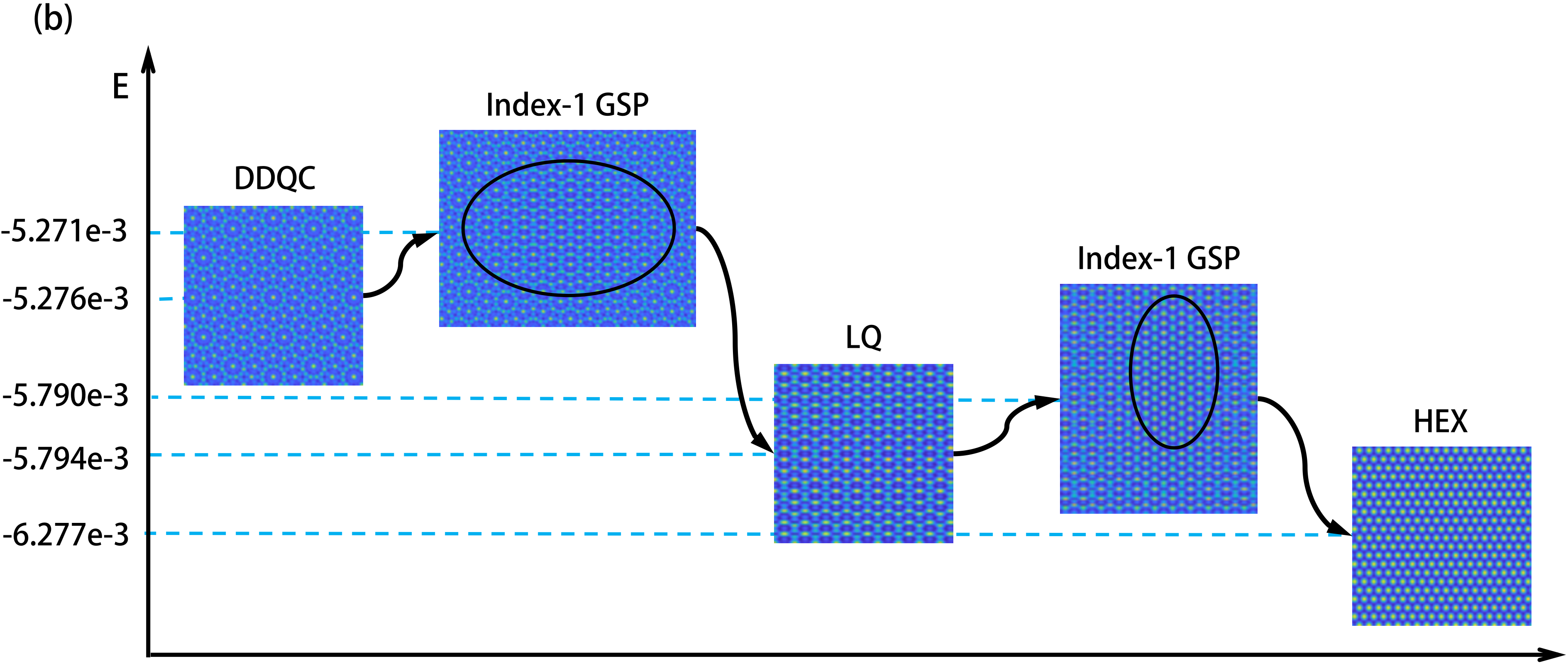}\\
        \end{minipage}%
        \caption{(a). Two-stage transition path from HEX to DDQC computed by the NPSS method in the LP model with $\varepsilon =5\times 10^{-6}$, $\alpha$ = $\sqrt{2}/2$. (b). Two-stage transition path from DDQC to HEX computed by the NPSS method in the LP model with $\varepsilon = 0.05$, $\alpha = 1$.}
        \label{fig:DDQCandHEX}
    \end{figure}

	Using the NPSS, we obtain
	a two-stage phase transition
	from HEX to DDQC, and an inverse phase transition from DDQC to HEX, as shown in \Cref{fig:DDQCandHEX}.
	These phase transitions have been also simulated in\,\cite{yin2021transition}. It is worth noting that multiple index-1 GSPs may be found starting from a GLM. The NPSS method could not guarantee to find the index-1 GSPs with the lowest energy barrier on the energy surface. However, if the targeted stable state is predicted, the NPSS method can be more efficient in locating the corresponding index-1 GSP with the target as a deterministic perturbation.
	
\begin{figure}[!hbpt]
		\centering
		\includegraphics[width=0.8\textwidth]{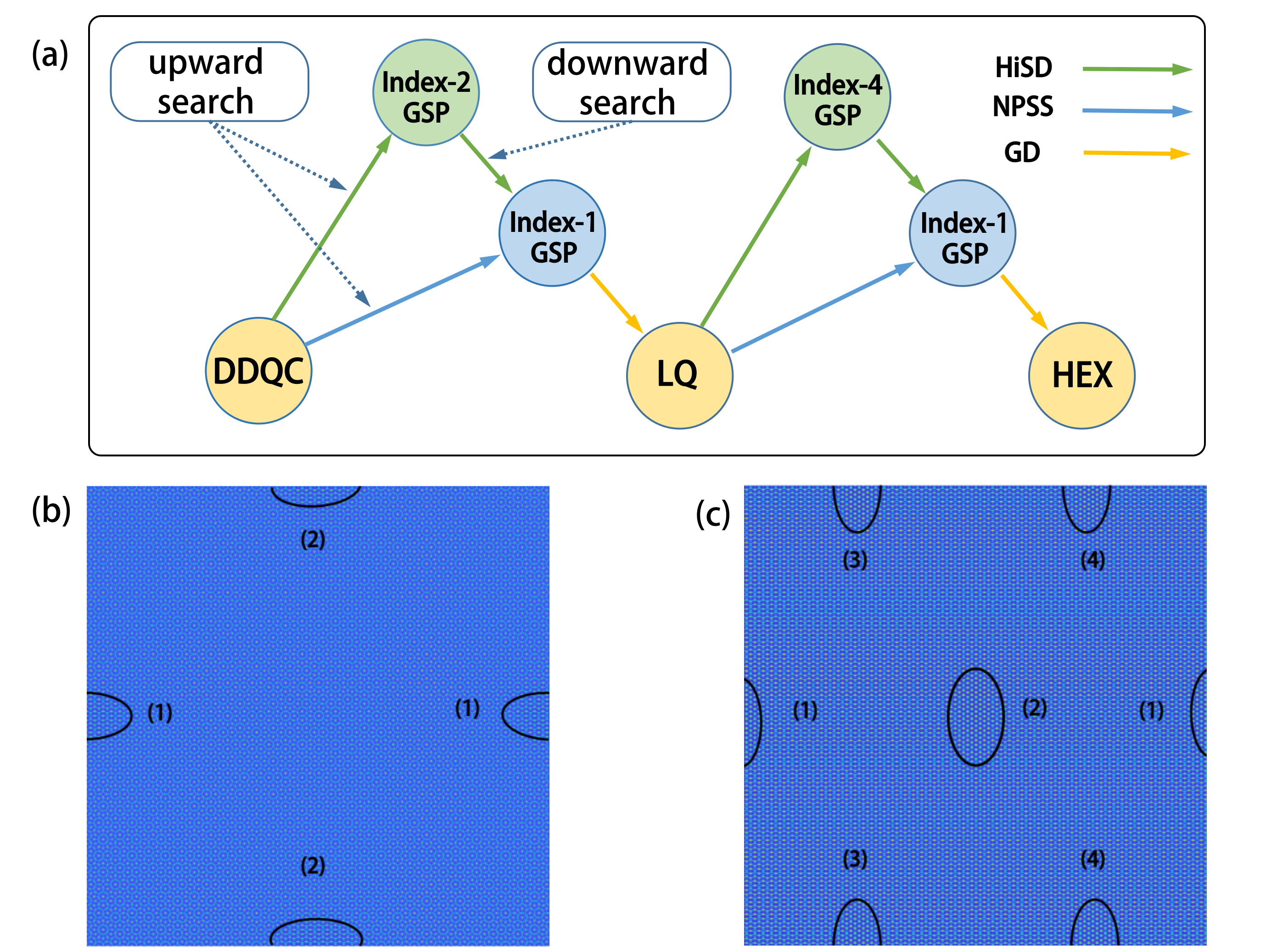}
		\caption{(a). Illustration of transition path searched by HiSD and NPSS method. (b). Index-2 GSP searched by HiSD method from DDQC. (c). Index-4 GSP searched by HiSD method from LQ. }
		\label{fig:illustration_DDQCtoHEX}
	\end{figure}

	Taking the phase transition from DDQC to HEX as an example, we use the HiSD method to search for the index-1 GSP from the same starting position $U^*$ with the same discrete format. \Cref{fig:illustration_DDQCtoHEX}\,(a) is a schematic diagram illustrating the phase transition process from DDQC to HEX, where HiSD and NPSS methods are applied to search for the index-1 GSP. Compared with the NPSS, which converges directly to the index-1 GSP in the upward search stage, the HiSD method first converges to the index-2 and index-4 GSPs (see \Cref{fig:illustration_DDQCtoHEX}\,(b) and (c)). It then requires a downward search to locate the index-1 GSP\,\cite{yin2021transition}.

	We compare the CPU times of both methods when searching for index-1 GSPs. As indicated in \Cref{tab:DDQCtoHEX_Efficiency}, the NPSS method is more efficient than the HiSD method when searching for index-1 GSPs. On average time in the upward search part, the HiSD method takes $3 \sim 5$ times longer than the NPSS method, depending on the dimension of the ascent space. Specifically, the HiSD method constructs and updates a $(l_k+1)$-dimensional ascent space $\mathcal{V} = {\rm span}\{\bm{v}^k_{1},\cdots,\bm{v}^k_{l_k},\bm{v}\}$, where $\{\bm{v}^k_{i} \}_{i=1}^{l_k}$ are eigenvectors of zero eigenvalues of $H(U)$, and $\bm{v}$ is a ascent direction. However, only the ascent direction $\bm{v}$ works for escaping from the basin. In contrast, the NPSS method provides a direct way to obtain $\bm{v}$, and only needs to update a 1-dimensional ascent space. Moreover, the HiSD method may converge to an index-$s$ GSP ($1 < s \leq l_k+1$), and must take more time to search down to the index-1 GSP.

	\begin{table}[!hbpt]
	    \centering
		\footnotesize{
		\caption{Comparison of the computational efficiency of the HiSD and NPSS methods when index-1 GSPs are searched from different initial states. US (DS) represents the upward (downward) search parts, respectively. Times, Iter, and Average time mean total CPU time, sum of iteration steps, and average CPU time per step, respectively.}
		\begin{tabular}{|c|c|c|c|c|c|c|}  
			\hline  
			Initial state          & Method  &  Convergence state &Times/s & Iter   & Average time/s      \\
			\hline  
			\multirow{3}{*}{DDQC}  & NPSS   & US: index-1 GSP &  4014    &  16240 & 0.25 \\ 
			\cline{2-6}
			                       & \multirow{2}{*}{HiSD}   & US: index-2 GSP &  21435   &  17230  & 1.24       \\
			
			                       &       &DS: index-1 GSP & 888 &  4370   & 0.20\\
			\hline
			\multirow{3}{*}{LQ}    & NPSS  & US: index-1 GSP  & 3334     &14460   & 0.23 \\  
			\cline{2-6}
				                   & \multirow{2}{*}{HiSD}   &US: index-4 GSP&   3805   &4030  & 0.70    \\
			                       &    & DS: index-1 GSP & 2643 & 8150    & 0.32\\	  
			\hline  
		\end{tabular}
		\label{tab:DDQCtoHEX_Efficiency}
		}
	\end{table}

	\section{Discussion}\label{sec:discussion}

\begin{figure}[H]
		\centering
		\includegraphics[width=0.9\textwidth]{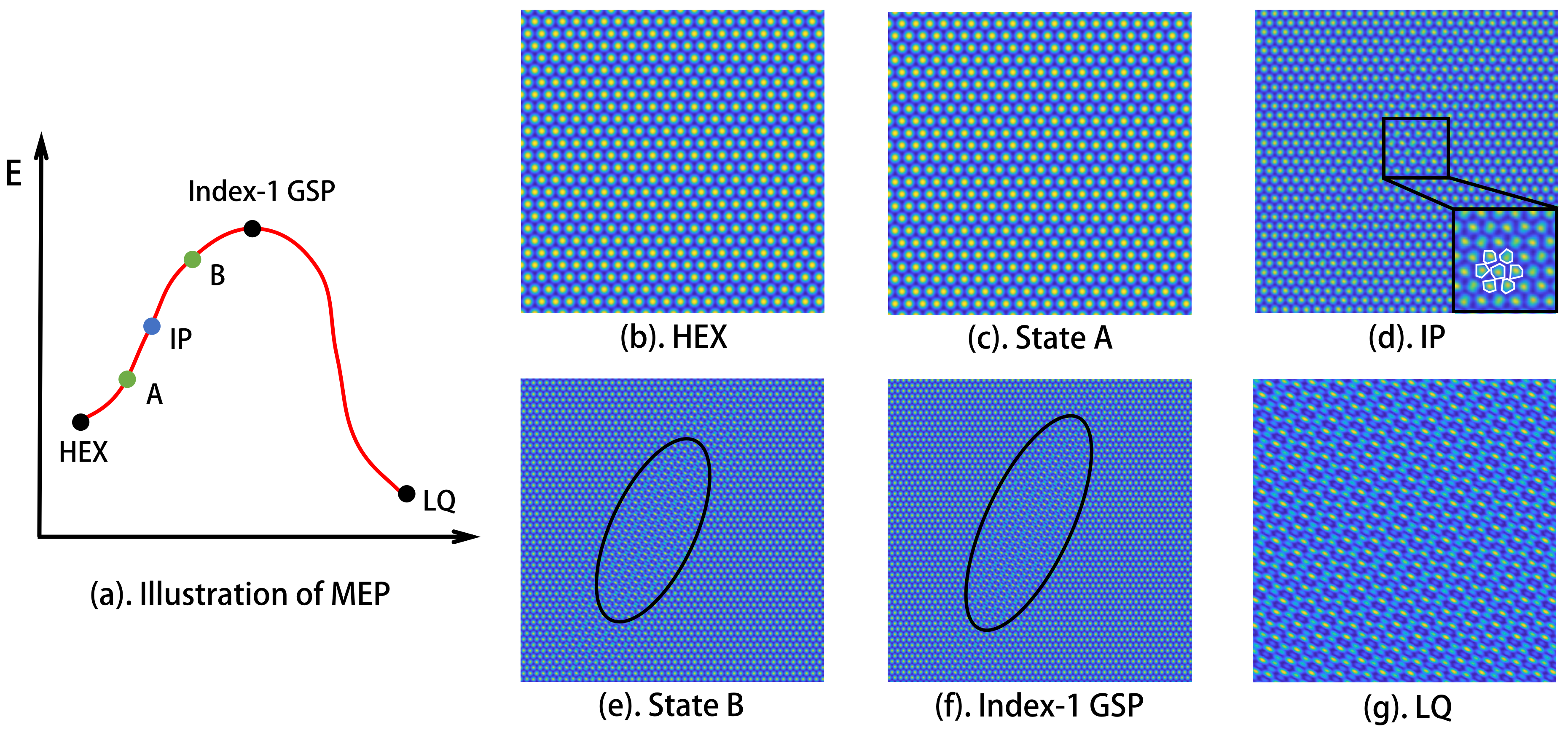}
		\caption{(a). Illustration of the MEP from HEX to LQ. (b)-(g). Several states on MEP marked in (a), where A and B are two intermediate states before and after reaching the IP, respectively. }
		\label{fig:states_on_MEP}
	\end{figure}

In this section, we investigate the relationship between symmetry-breaking and nullspaces during the phase transition. Considering the phase transition from HEX to LQ in \Cref{fig:DDQCandHEX}\,(a), we study the structural changes on MEP. There is an IP on MEP where the smallest positive eigenvalue becomes negative. As shown in \Cref{fig:states_on_MEP}, before reaching the IP, the structure like state A on the MEP is not obviously different from that of HEX. At the IP, the cylinders transform into a polygonal shape as shown in (d). After passing the IP, a new cluster gradually appears until the critical nucleus is formed, as shown in (e) and (f), respectively. It can be seen that the IP is the critical point where the symmetry-breaking begins to occur. \Cref{sec:periodicity} has shown that the translational invariance of the stable state $U$ is related to the dimension of the nullspace $\mathcal{W}^k(U)$. 
Does the invariance of the ordered structure before the IP on MEP also correspond mathematically to the invariance of the nullspace? A natural question is whether the states on the MEP before the IP can maintain the nullspace invariance.

To further explore the property of nullspace invariance on the MEP, we introduce some notations and definitions below. We generate a sequence of states $\{U^{(j)}\}_{j=0}^n$ on MEP from GLM to index-1 GSP by the GD method, where $U^{(0)}$, $U^{(m)}$ and $U^{(n)}$ ($0 < m < n$) represent the GLM, IP and index-1 GSP, respectively. We then compute the smallest six eigenvalues of $H(U^{(j)}), j= 1,2,\cdots,n$, and their corresponding eigenvectors by LOBPCG.
We rewrite $\mathcal{W}^{k}(U^{(0)}) = \text{span}\{\bm{v}^k_{0,1}, \cdots, \bm{v}^k_{0,l_k}\}$. $\bm{v}$ is one of eigenvectors on the MEP. We define an angle $\theta$ to measure relations between $\bm{v}$ and $\mathcal{W}^{k}(U^{(0)})$ as
	\begin{align}
	    \begin{aligned}
		\theta &= \arccos(\|\bm{v_{p}}\|_{\ell^2}/\|\bm{v}\|_{\ell^2}),\\
		\bm{v_{p}} &= \sum^{i=l_k}_{i=1}\frac{\langle \bm{v}, \bm{v}^k_{0,i}\rangle}{\langle \bm{v}^k_{0,i}, \bm{v}^k_{0,i} \rangle} \bm{v}^k_{0,i},\notag
	\end{aligned}
	\end{align}
	where $\bm{v_{p}}$ is the projection of $\bm{v}$ in $\mathcal{W}^{k}(U^{(0)})$.
	
\begin{figure}[H]
		\centering
		\includegraphics[width=1\textwidth]{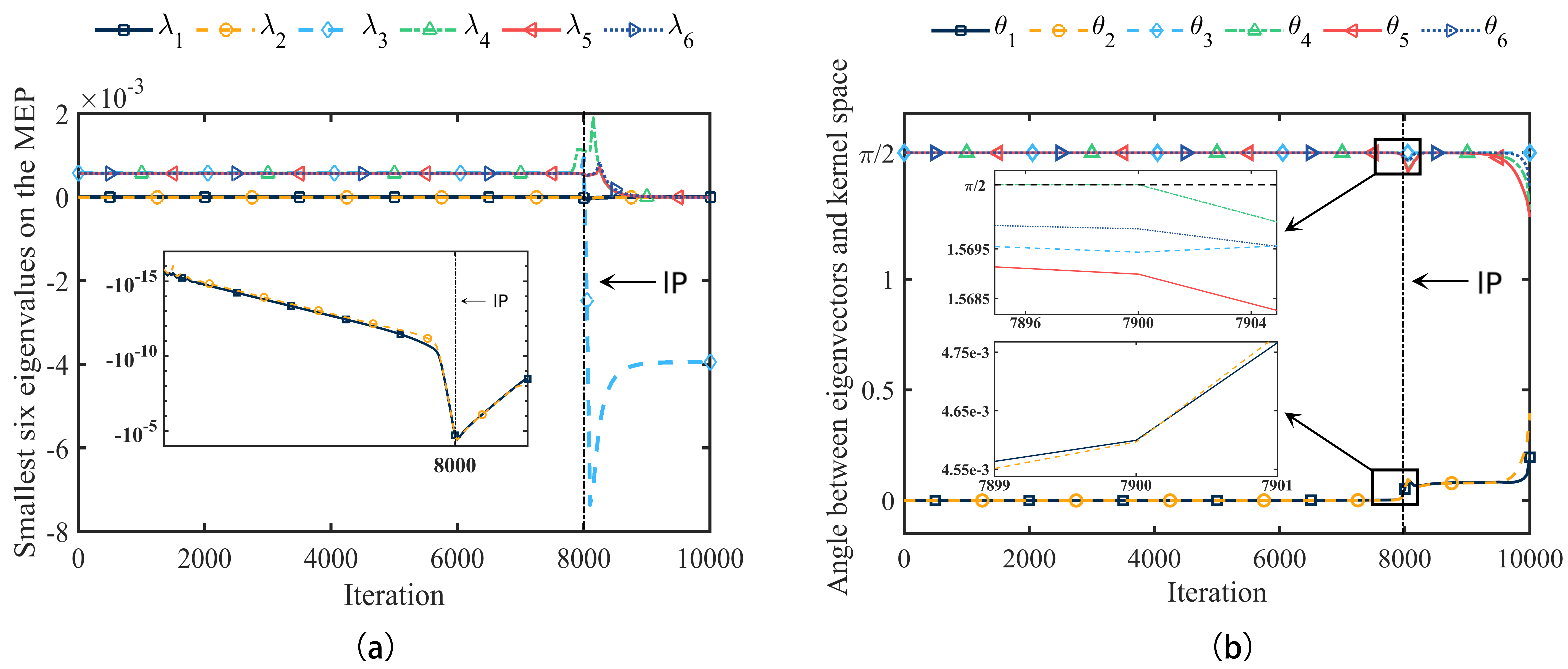}
		\caption{(a). Variation of the smallest six eigenvalues on the MEP from HEX to the index-1 GSP. The subgraph shows the variation details of $\lambda_1$ and $\lambda_2$. (b). Angles $\{\theta_i\}^6_{i=1}$ between the nullspace of HEX and the eigenvectors corresponding to the smallest six eigenvalues on the MEP. The subgraphs show the details of angles near the IP. }
		\label{fig:eigenvalues_and_angles_HEXtoLQ}
	\end{figure}

	Then we study properties of nullspaces $\{\mathcal{W}^k(U^{(j)})\}^{n}_{j=0}$ by exploring smallest six eigenvalues $\{\lambda_i\}_{i=1}^6$ and angles $\{\theta_i\}^{6}_{i=1}$ between corresponding eigenvectors and $\mathcal{W}^k(U^{(0)})$.    \Cref{fig:eigenvalues_and_angles_HEXtoLQ}\,(a) illustrates the variation curves of the six smallest eigenvalues from the GLM to the index-1 GSP. From $U^{(0)}$ to the neighborhood of the IP, $|\lambda_1|$ and $|\lambda_2|$ remain below $\epsilon_0$ = 1e-10, which can be considered as zero eigenvalues. This means that the dimension of nullspaces along the MEP before reaching the neighborhood of the IP is equal to that of the GLM.     
	As shown in \Cref{fig:eigenvalues_and_angles_HEXtoLQ}\,(b), $\theta_1$ and $\theta_2$ are less than 1e-2 before reaching the neighbourhood of IP, indicating that the nullspace on the MEP belongs to $\mathcal{W}^{k}(U^{(0)})$ in this part of the MEP.
	In summary, the nullspace of the state on the MEP from the GLM to the neighborhood of the IP is approximately equal to the nullspace of the GLM, consistent with the structural invariance of the parent phase during the phase transition from the GLM to the IP.
	
    We refer to this feature as the nullsapce-preserving property, and surmise that symmetry preservation implies nullspace preservation in phase transitions with translational invariance. The IP on the MEP corresponds mathematically to the dividing point where the nullspace is no longer preserved, and physically to the critical point where the symmetry begins to change. This property implies that the NPSS method may not need to segment and update nullspaces multiple times during escaping from the basin for studying phase transition with translational invariance. Computational efficiency may be able to be further developed.

%	\newpage
	\section{Conclusion}\label{sec:conclusion} 
	In this paper, we develop the NPSS method to study phase transitions with translational invariance based on the LB and LP models. The NPSS method can quickly escape from the basin and avoid the effects of degenerate problems. The NPSS method divides the process of escaping from the basin into multiple segments and climbs upward along an ascent direction orthogonal to the nullspace of an initial state in each segment. This method can reduce the computational costs for updating the nullspace at every step and ensures the ascent direction is efficient. We employ the NPSS method to study the nucleation and phase transition between quasicrystal and crystal, and between crystals. Results show that the NPSS method can locate index-1 GSPs with an economical computational cost compared with HiSD, and it can efficiently reveal transition paths between ordered structures with translational invariance. In particular, we identify an important critical point, the IP, on the MEP associated with symmetry-breaking. Before reaching the IP, we find that the nullspaces on the MEP have the $nullspace$-$preserving$ property, which may be able to promote further development of NPSS methods in computational efficiency. In the future, we plan to apply the NPSS method to more interesting phase transition problems and to improve it with problem-dependent features.

    \section*{Acknowledgement} 
	We would like to thank Prof.~Lei Zhang and Dr.~Jianyuan Yin for their valuable suggestions and for sharing their HiSD algorithm's codes. 
    We also would like to thank the anonymous reviewers and the associate editor for their constructive comments in improving this work.
	
	\bibliographystyle{siamplain}
	\bibliography{references}

\end{document}